%% file: ss_spde_v3.tex
\documentclass[12pt]{article}
\usepackage{natbib,bm}
 \usepackage{amstext}
\author{Yu Ryan Yue}
\usepackage{graphicx}
\usepackage{amsmath}
\usepackage{pdflscape}
\usepackage{subfigure}

\input Macros

\begin{document}

\def\spacingset#1{\renewcommand{\baselinestretch}%
{#1}\small\normalsize} \spacingset{1}

  \title{\bf Bayesian Adaptive Smoothing Spline using Stochastic Differential Equations}

\author{Yu Ryan Yue$^{\rm a}$\thanks{Corresponding author. Email: yu.yue@baruch.cuny.edu
}, ~Daniel Simpson$^{\rm b}$, Finn Lindgren$^{\rm c}$ \& H\aa vard Rue$^{\rm b}$\\\\
$^{\rm a}${\em{Baruch College, The City University of New York, USA.}}\\ 
$^{\rm b}${\em{Norwegian University of Science and Technology, Norway.}}\\
$^{\rm c}${\em{University of Bath, UK.}}\\
}

  \maketitle

\bigskip
\begin{abstract}
The smoothing spline is one of the most popular curve-fitting methods, partly because of empirical evidence supporting its effectiveness and partly because of its elegant mathematical formulation. However, there are two obstacles that restrict the use of smoothing spline in practical statistical work. 
Firstly, it becomes computationally prohibitive for large data sets because the number of basis functions roughly equals the sample size. Secondly, its
global smoothing parameter can only provide constant amount of smoothing, which often results in poor performances when estimating inhomogeneous 
functions. In this work, we introduce a class of adaptive smoothing spline models that is derived by solving certain stochastic differential equations with finite element methods. The solution extends the smoothing parameter to a continuous data-driven function, which is able to capture the change of the smoothness of underlying process. The new model is Markovian, which makes Bayesian computation fast. A simulation study and real data example are presented to demonstrate the effectiveness of our method. 
\end{abstract}

\noindent%
\spacingset{1.45}

\newpage

\section{Introduction}\label{intro}
The smoothing spline is one of the most popular nonparametric regression methods, partly because of empirical evidence supporting its effectiveness and partly because of its elegant mathematical formulation. Consider the model
\be\label{mod:nonpar}
y_i = f(t_i) + \vareps_i,\quad i=1,\ldots,n; \quad t_i\in\mathcal{T},
\ee
where $\bfy=(y_1,y_2,\ldots,y_n)$ is the vector of observations, $f$ is some ``smooth" function defined on some index set $\mathcal{T}$, and $\vareps_i\stackrel{iid}{\sim}N(0,\tau^{-1})$ with precision (inverse of variance) $\tau$. The smoothing spline of degree $2p-1$ is defined as the solution of the following minimization problem,
\be\label{min:ss}
\hat{f}=\mbox{arg}\min_{f}\bigg[\sum_{i=1}^{n}\Big(y_{i}
-f(t_{i})\Big)^{2}+\lambda\int_{\mathcal{T}}\Big(f^{(p)}(t)\Big)^{2}dt\bigg],
\ee
where $\lambda>0$ is the {\em smoothing parameter} and $ f^{(p)}(t)$ is the $p$th derivative of $f(t)$. The parameter $\lambda$ controls the trade-off between fidelity to the data in terms of the residual sum of squares against smoothness of the fit in terms of the integrated squared derivative. The value of $p$ is often taken to be 1 or 2, corresponding to linear and cubic smoothing spline, respectively. From frequentist point of view, the solution $\hat f$ can be explicitly derived within a reproducing kernel Hilbert space and $\lambda$ is usually estimated via cross-validation or generalized cross-validation method \citep[see e.g.,][]{Wahb:spli:1990,Gu:smoo:2002}. From Bayesian point of view, the $\hat f$ is the mean of the posterior distribution of $f$ yielded by taking a partially improper Gaussian prior taken on the function space \citep{wahba78,Euba:nonp:1999,Spec:Sun:full:2003}.



There are two obstacles that restrict using smoothing spline estimators in practical statistical work. Firstly, they become computationally intractable for large data sets because the number of basis functions roughly equals the sample size \citep{Wahb:spli:1990, Gree:Silv:nonp:1994}. The second obstacle stems from the smoothing parameter $\lambda$. A single parameter $\lambda$ implies that the underlying mean process $f(t)$ has a constant amount of smoothing, which is not always realistic in practice. It often results in the poor performance of smoothing spline, especially when estimating inhomogeneous functions.

To overcome the computation issue, one approach is to employ regression splines \citep[see][for a comprehensive review]{hansen02}. The basis implied by solving the spline smoothing problem for a small representative data set is found and this small basis is used to construct a model for the full data set of interest. The model is typically fitted as a linear or generalized linear model without imposing a roughness penalty. The covariate points that are used to obtain the reduced basis are known as the `knots' of the regression spline. The number of knots controls the flexibility of the model, but unfortunately their locations tend to have a marked effect on the fitted model. Some of the problems with knot placement can be partially alleviated by using penalized regression splines (P-splines), where the required penalty is associated with the regression spline basis. It is interesting to note that there are two versions of P-splines, which can be distinguished by their bases and penalties in use. \cite{Eile:Marx:pspline:1996} introduced the P-splines with B-spline basis and differencing penalty, while \cite{Rupp:Carr:spat:2000} and \cite{rupp:wand:carr:2003} proposed a competing method with truncated power basis and ridge penalty. Both P-splines have recently gained incredible popularity in statistics and applied fields due to their easy implementation using linear mixed model formulation \citep{Eilers10}. \cite{OSullivan86} also introduced a similar penalized spline approach using B-spline basis, but with a more complicated penalty derived from the integrated squared derivative of the fitted curve. The O'Sullivan spline was recently revived by \cite{wand08}, who showed that it possess attractive features, e.g., smoothness, numerical stability and natural boundary properties. \cite{simp:helt:lind:12} characterized
the connection between O'Sullivan splines, classical smoothing splines
and the Markovian models considered in this paper.


To increase its smoothing flexibility, many authors have proposed to make smoothing splines adaptive, e.g., local generalized cross-validation approach in \cite{Cumm:Fill:Nych:2001}, adaptive $L$-splines in \cite{Abra:Stei:impr:1996}, hybrid adaptive splines in \cite{luo:wahba:97}, and spatially adaptive smoothing splines in \cite{Pin:Spec:Holm:spat:2006}. There is also extensive literature on adaptive P-splines, where a functional structure on the smoothing parameters is imposed in the ordinary P-spline models. The adaptive smoothing function is often chosen as another layer of P-spline with a set of subknots. Typical works include \cite{Lang:Brez:baye:2004}, \cite{Bala:Mall:Carr:spat:2005}, \cite{brezger:csda:06}, \cite{Crai:Rupp:07}, \cite{Kriv:Crai:Kaue:fast:2008} and \cite{scheipl:csda:09}. As their ordinary counterparts, the adaptive P-splines need ``good'' knots and subknots to provide appropriate adaptive smoothing. Several other spline-based adaptive smoothing methods are proposed as well, including local polynomial models with adaptive window widths \citep{fan96}, adaptive regression splines \citep{denison98,zhou01,Dim:Gen:Kas:baye:2001,Holm:Mall:bay:2001} and mixtures of smoothing splines \citep{Wood:Jiang:Tan:mix:2002,wood:kohn:08}

In this work, we propose a unified and efficient Bayesian approach to model smoothing splines, which can be easily equipped with adaptive smoothing feature. The method is based on constructing Gaussian Markov random field (GMRF) representations for adaptive smoothing splines by solving certain stochastic differential equations. We here provide a brief introduction for GMRF.  A random vector $\bfw=(w_1,\ldots,w_n)'$ is a GMRF if it has density of form 
\be\label{gmrf_org}
[\bfw\mid\delta] \propto |\delta\bfQ|^{1/2}_+\exp\left(-\fr{\delta}{2}(\bfw-\bfmu)'\bfQ(\bfw-\bfmu)\right),
\ee 
where $\delta>0$ is scale parameter, $\bfmu$ is mean vector, and $\bfQ$ is so-called {\em precision} matrix. The notation $|\bfA|_+$ denotes the generalized determinant of matrix $\bfA$, which is the product of its nonzero eigenvalues. The full conditionals $\pi(w_i\mid\bfw_{-i})$, $i=1,\ldots,n$, only depend on a set of neighbors $\mathcal{N}_i$ to each site $i$. The computational gain comes from the fact that the zero-pattern of matrix $\bfQ$ relates directly to the notion of neighbors: $Q_{ij}\ne 0 $ if and only if $i\in\mathcal{N}_j\cup j$ \citep[see e.g.,][Sec 2.2]{GMRFbook}. The GMRFs allow for fast direct numerical algorithms, as numerical factorization of $\bfQ$ can be done using sparse matrix algorithms at a typical cost of $\mathcal{O}(n)$; see \cite{GMRFbook} for detailed algorithms. Such good computational properties are of major importance in Bayesian inferential methods. This is further enhanced by the link to nested integrated Laplace approximations (INLA) \citep{Rue:Mart:Chop:inla:2009}, which allows for fast and accurate Bayesian inference for latent Gaussian field models.

The connection between GMRF and smoothing splines have been explored by several authors. \cite{Spec:Sun:full:2003} showed that the random walk (RW) models (a subclass of GMRF) \citep[e.g.,][]{Fahr:Wage:smoo:1996, fahr:knor:00, Fahr:Lang:baye:2001}, can be used as priors to derive the discretized Bayesian smoothing spline estimator. \cite{Lang:Fron:Fahr:func:2002} and \cite{yue:spec:sun:aism:12} made the RW models spatially adaptive by introducing local smoothing parameters into the models. However, all the RW models mentioned above are only appropriate for the data observed at regular locations. \cite{lind:rue:08} considered a second-order RW (RW2) model as a discretely observed continuous time process, which is derived by solving a stochastic differential equation (SDE) with finite element method. The resulting RW2 model is resolution consistent and has a GMRF representation of the cubic smoothing spline, with equally good performance but more computational efficiency.



The aim of this paper is to extend \citeauthor{lind:rue:08}'s work in regard to spatial adaptation. More specifically, we enable their RW2 model to be spatially adaptive by carefully adding a smoothing function to the SDE. The smoothing function is able to provide various amounts of smoothing as required by the data. The solution of this modified SDE is thus a spatially adaptive smoothing spline, whose GMRF representation is explicitly available for any collection of locations. Compared to the existing methods, the adaptive smoothing models considered in this paper have a number of advantages. In particular, they have both a convenient computational
form and a well-understood continuous limit. This not only allows for fast
computation, but also provides the comfort
that issues like knot spacing will only have a minimal and well-known effect on the model \citep[see][for a discussion]{simp:helt:lind:12}. Furthermore, they provide a satisfactory extension of the models
in \cite{lind:rue:08} to adaptive smoothing, which means that we
can use the intuition built off those models, and correspondingly off
RW2 models on regularly-spaced knots, to understand these models.

\section{Bayesian smoothing spline using SDE}\label{sec:bss}



\cite{kim:wah:70} and \cite{wahba78} showed that the smoothing spline $\hat f$ in (\ref{min:ss}) is equivalent to Bayesian estimation with a partially improper prior generated by the following stochastic differential equation (SDE)
\be\label{sde:nad}
d^pf(t)/dt^p = dW(t)/dt,
\ee
where the function $W(t)$ is a zero mean Wiener process with variance $t$, and $dW(t)/dt$ is often referred to as ``white noise". Letting $(t)_+=t$ for $t\ge0$ and $(t)_+ = 0$ otherwise, the exact solution of SDE (\ref{sde:nad}) is shown to be
\be\label{soln:sde}
f(t) = \beta_0 + \beta_1t + \cdots + \beta_{p-1}t^{p-1}+ Z(t)/\sqrt{\delta},\quad t\in\mathcal{T},
\ee
where $\delta>0$, $\beta_0,\beta_1,\ldots,\beta_p\sim N(0, \xi)$ as $\xi\goto\infty$, and $Z(t)$ is a zero mean Gaussian stochastic process with $E[Z(s)Z(t)] = \Sigma(s,t)$ and
\ba
\Sigma(s,t) = \int_0^1\fr{(s-u)^{p-1}_+}{(p-1)!}\fr{(t-u)_+^{p-1}}{(p-1)!}du.
\ea
We actually take a partially improper prior on $f$, which is ``diffuse" on the coefficients of the polynomials of degree $p-1$, and ``proper" over the random process $Z(t)$. Then, the $\hat f$ has the property $\hat f(t) = \lim_{\xi\goto\infty}E_\xi\{f(t)\mid\bfy,\tau,\delta\}$, which is the expectation over the posterior distribution of $f(t)$ with the prior defined in (\ref{soln:sde}). Note that the smoothing parameter $\lambda$ now becomes $\lambda = \delta/\tau$. After taking sensible priors on $\tau$ and $\delta$, the fully Bayesian inference on $\hat f$ can be straightforwardly carried by Monte Carlo Markov chain (MCMC) method \citep{Spec:Sun:full:2003,yue:spec:sun:aism:12}


Unfortunately, the prior (\ref{soln:sde}) is computationally intensive for large data sets because the covariance matrix of $Z(t)$ is completely dense. We therefore solve SDE (\ref{sde:nad}) using a finite element approach as introduced in \cite{lind:rue:08}. The solution will be shown to be a GMRF of form in (\ref{gmrf_org}). Note that we here only consider cubic smoothing spline ($p=2$), which is well known to provide the best overall performance. Let $t_1<t_2<\cdots<t_n$ be the set of fixed points, which are often observed locations, but do not have to be. Define the inner product
$
\langle f, g \rangle= \int f(t)g(t) dt
$, 
where the integral is over the region of interest. We seek a stochastic weak solution of (\ref{sde:nad}) for $p=2$ that satisfies
\be\label{sde:weak}
\left\langle \phi, d^2f/dt^2\right\rangle \stackrel{d}{=} \left\langle \phi, dW/dt\right\rangle
\ee
for {\em any} sensible test function $\phi(t)$, where $\stackrel{d}{=}$ denotes equality in distribution \citep{walsh86}.  It is impossible to test (\ref{sde:weak}) against every function $\phi(t)$, so we chose a finite set $\{\phi_i(t)\}_{i=1}^n$ instead.

We then construct a finite element representation of $f(t)$ as
\be\label{approx:f}
f(t) \approx \sum_{j=1}^n \psi_j(t) w_j,
\ee
for some chosen basis functions $\psi_j$ and random weights $w_j$. Letting $h_j = t_{j+1}-t_j$ for $j=1,\ldots,n-1$, a common choice of basis is the piecewise linear functions
\ba
\psi_j(t) = \left\{
\begin{array}{ll}
0, & t < t_{j-1},\\
\fr{1}{h_{j-1}}(t - t_{j-1}), & t_{j-1}\le t< t_j,\\
1 - \fr{1}{h_j}(t-t_j), & t_j\le t < t_{j+1},\\
0, & t_{j+1}\le t.\\
\end{array}
\right.
\ea
An interpretation of the representation (\ref{approx:f}) with this chosen basis functions is that the weights determine the values of the field at the locations, and the values in the interior of the intervals are determined by linear interpolation. The full distribution of the continuously indexed solution is determined by the joint distribution of the weights $\bfw=(w_1,\ldots,w_n)^T$. 

%
Finally, we let the test functions be the same as our basis functions, which is known as {\em Galerkin} finite element method. Substituting (\ref{approx:f}) into (\ref{sde:weak}) for this set of test functions, we end up with a system of linear equations
\be\label{sde:weak:finite}
\sum_{j=1}^nw_j\left\langle \psi_i, d^2\psi_j/dt^2\right\rangle\stackrel{d}{=}  \big\langle \psi_i,  dW/dt\big\rangle,\quad i=1,\ldots,n.
\ee
The finite dimensional solution is obtained by finding the distribution of $\bfw$ that fulfills the weak SDE formulation (\ref{sde:weak:finite}). It can be shown that the left hand side of (\ref{sde:weak:finite}) can be written as $\bfH\bfw$,
where $\bfH$ is an $n\times n$ tridiagonal matrix whose non-zero entries are 
\be\label{sde:lhs}
\bfH[i,i-1]= \fr{1}{h_{i-1}}, ~\bfH[i,i] = -\left(\fr{1}{h_{i-1}} + \fr{1}{h_i}\right), ~\bfH[i,i+1]= \fr{1}{h_i}
\ee
for $2\le i \le n-1$, 
since $\psi_i$ only overlap for neighboring basis functions. The entries of the first and last row in $\bfH$ are zeroes. Given the statistical properties of white noise, the inner product on the right-hand side of (\ref{sde:weak:finite}) is a Gaussian distribution with zero mean and covariance matrix $\bfB=[\langle\psi_i,\psi_j\rangle]_{i,j=1}^n$, whose nonzero entries are given by
\ba
\bfB[i,i-1]=\fr{h_{i-1}}{6},~~\bfB[i,i] = \fr{h_{i-1}+h_i}{3},~~\bfB[i,i+1] = \fr{h_i}{6},
\ea 
with modifications at the boundaries. To achieve distribution equality in (\ref{sde:weak:finite}), the random vector $\bfw$ has the density of form (\ref{gmrf_org}) with $\bfmu = \bfzero$ and $\bfQ = \bfH'\bfB^{-1}\bfH$. However, such $\bfQ$ is the dense matrix due to the dense $\bfB^{-1}$, making the Galerkin model computationally expensive. Lindgren and Rue (2008) showed that without changing the solution we may replace $\bfB$ by a diagonal matrix $\tilde{\bfB}$ with $\tilde{\bfB}[i,i]=\langle \psi_i, 1\rangle$, giving 
\be\label{mat:tildeB}
\tilde{\bfB}[1,1] = \fr{h_1}{2}, ~~\tilde{\bfB}[i,i]=\fr{h_{i-1} + h_i}{2}, ~~\tilde{\bfB}[n,n] = \fr{h_{n-1}}{2}. 
\ee
As a result, the matrix $\bfQ = \bfH'\tilde\bfB^{-1}\bfH$ becomes sparse and $\bfw$ is thus a GMRF. 
It is straightforward to verify that $\bfQ$ has rank $n-2$, with the null space spanned by vectors $(1,\ldots,1)^T$ and $(t_1,\ldots,t_n)^T$. It indicates that the resulting field is invariant to addition of a linear trend, coinciding with the result obtained by \cite{wahba78} for cubic smoothing spline. 

We have now derived a GMRF $\bfw$ as the weights of a basis function expansion (\ref{approx:f}), which approximates the continuous function $f(t)$ everywhere. \cite{simp:helt:lind:12} showed that the convergence of the approximation depends solely on the basis functions. Given any set of enough points $t_i$, using the piecewise linear functions yields the best finite approximation to the continuos process regardless of their locations. Also, the method described above works for any set of test and basis functions when all of the computations make sense. Actually, \citeauthor{simp:helt:lind:12} showed that the O'Sullivan spline can be exactly derived by solving the SDE in (\ref{sde:nad}) using cubic B-splines as basis functions and their second derivatives as test functions. However, one should be aware that the wrong choice of global basis functions will destroy the Markov structure, and not all sets of basis functions will provide good approximations to $f(t)$. 

\section{Extensions to adaptive smoothing spline}\label{sec:bass}
Besides their intriguing theoretical and computational properties, one of the most exciting aspects of the SDE spline models is their flexibility: it is straightforward to extend them to adaptive smoothing spline models. The basic idea is that by making the smoothing parameter vary in space, we will be able to control the local
smoothing properties of the spline.
We here present two different adaptive SDE formulations, from both of which we are able to derive the GMRF models that provide appropriate adaptive smoothing.

\subsection{Adaptive SDE I}
One way to extend SDE (\ref{sde:nad}) is as follows:
\be\label{sde}
\lambda(t)d^2f(t)/dt^2 =  dW(t)/dt,
\ee
where the positive $\lambda(t)$ can be seen as an adaptive smoothing function, compared to the global smoothing parameter $\lambda$ in ordinary smoothing splines. A small $\lambda(t)$ allows big second derivative of $f(t)$ for roughness, while a large value diminishes the derivative to increase smoothness. The solution to (\ref{sde}) is related to the spatially adaptive smoothing spline introduced in \cite{Pin:Spec:Holm:spat:2006}, minimizing
\be\label{min:adss}
\sum_{i=1}^n\Big(y_i - f(t_i)\Big)^2 + \int\left[\lambda(t) f''(t)\right]^2 dt.
\ee
Using a piecewise-constant model for $\lambda(t)$, \citeauthor{Pin:Spec:Holm:spat:2006}\ derived closed-form solutions for the corresponding reproducing kernels of the Hilbert space. Their method, however, is computationally intensive since the matrix of reproducing kernel is completely dense.

Following the non-adaptive case, we seek a weak solution of (\ref{sde}) by achieving 
\be\label{sde:ad:weak}
\big\langle \psi_i, \lambda d^2f/dt^2\big\rangle \stackrel{d}{=}\big\langle \psi_i, dW/dt\big\rangle, \quad i=1,\ldots,n.
\ee
Using the basis representation in (\ref{approx:f}) as well as Galerkin approximation, the left hand side of (\ref{sde:ad:weak}) can be proved to be $\bfLambda\bfH\bfw$,
where $\bfLambda$ is a diagonal matrix of $\bflambda = (\lambda(t_1),\ldots,\lambda(t_n))^T$ and $\bfH$ is the matrix as in (\ref{sde:lhs}) (see Appendix for the proof). Since the right-hand side of (\ref{sde:ad:weak}) is the same as in (\ref{sde:weak:finite}), the $\bfw$ is also a GMRF with zero mean and the following precision matrix
$$
\bfQ_\lambda = \bfH'\bfLambda\tilde\bfB^{-1}\bfLambda\bfH.
$$
It is easy to see that $\bfQ_\lambda$ is symmetric and banded with non-zero entries of $i$th row given by
\ba
&&\bfQ_\lambda[i,i-2] = \fr{2\lambda^2(t_{i-1})}{h_{i-2}h_{i-1}(h_{i-2} + h_{i-1})}, ~
\bfQ_\lambda[i,i-1] = -\fr{2}{h^2_{i-1}}\left(\fr{\lambda^2(t_{i-1})}{h_{i-2}} + \fr{\lambda^2(t_i)}{h_i}\right),\\
&& \bfQ_\lambda[i,i] = \fr{2\lambda^2(t_{i-1})}{h^2_{i-1}(h_{i-2} + h_{i-1})} + \fr{2\lambda^2(t_i)}{h_{i-1} h_i}\left(\fr{1}{h_{i-1}}+\fr{1}{h_i}\right) + \fr{2\lambda^2(t_{i+1})}{h_i^2(h_i + h_{i+1})}.
\ea
At the discretization boundaries, we use the convention that terms with non-existing components are ignored, that is $h_{-1} =h_0 =h_n =h_{n+1} =\infty$. This affects only the upper left and lower right corner of $\bfQ_\lambda$ as follows:
\ba
&& \bfQ_\lambda[1,1]=\fr{2\lambda^2(t_2)}{h_1^2(h_1+h_2)}, ~~\bfQ_\lambda[2,1]= -\fr{2\lambda^2(t_2)}{h_1^2h_2},\\
 &&\bfQ_\lambda[2,2]=\fr{2\lambda^2(t_3)}{h_2^2(h_2+h_3)} + \fr{2\lambda^2(t_2)}{h_1h_2}\left(\fr{1}{h_1} + \fr{1}{h_2}\right),\\
&& \bfQ_\lambda[n-1,n-1]=\fr{2\lambda^2(t_{n-2})}{h_{n-2}^2(h_{n-3}+h_{n-2})} + \fr{2\lambda^2(t_{n-1})}{h_{n-2}h_{n-1}}\left(\fr{1}{h_{n-2}} + \fr{1}{h_{n-1}}\right),\\
&&\bfQ_\lambda[n, n]=\fr{2\lambda^2(t_{n-1})}{h_{n-1}^2(h_{n-2}+h_{n-1})}, ~~\bfQ_\lambda[n,n-1] = -\fr{2\lambda^2(t_{n-1})}{h_{n-2}h_{n-1}^2}.
\ea
Note that $\bfQ_\lambda$ does not involve $\lambda(t_1)$ or $\lambda(t_n)$ because the first and last rows of $\bfH$ are zeroes.


\subsection{Adaptive SDE II}
An alternative SDE that we can use for adaptive smoothing is 
\be\label{sde2}
d^2\lambda(t)f(t)/dt^2 = dW(t)/dt,
\ee
where $\lambda(t)$ can be seen as a instantaneous variance or local scaling, which compress and stretch the function. A small $\lambda(t)$ compresses the scale giving quick oscillations, while a high value stretch $f(t)$, decreasing the roughness. Adopting notation $\tilde f(t) = \lambda(t)f(t)$, formulation (\ref{sde2}) corresponds to minimizing
\ba
\sum_{i=1}^n\Big(y_i - \tilde f(t_i)\Big)^2 + \int\tilde f''(t)^2 dt.
\ea

The weak solution of (\ref{sde2}) can also be found using Galerkin method to satisfy
\be\label{sde:weak2}
\big\langle \psi_i,  d^2\lambda f/dt^2\big\rangle \stackrel{d}{=}\big\langle \psi_i, dW/dt\big\rangle, \quad i=1,\ldots,n,
\ee
whose left-hand side can be written as
$\bfH\bfLambda\bfw$, where $\bfH$ and $\bfLambda$ are defined as above (see Appendix for the proof). Again, the $\bfw$ is a GMRF with zero mean and precision matrix 
\ba
\bfQ_\lambda = \bfLambda\bfH'\tilde\bfB^{-1}\bfH\bfLambda,
\ea
whose nonzero entries can be explicitly written out as
\ba
&&\bfQ_\lambda[i,i-2] = \fr{2\lambda(t_{i-2})\lambda(t_i)}{h_{i-2}h_{i-1}(h_{i-2} + h_{i-1})}, ~
\bfQ_\lambda[i,i-1] = -\fr{2\lambda(t_{i-1})\lambda(t_i)}{h^2_{i-1}}\left(\fr{1}{h_{i-2}} + \fr{1}{h_i}\right),\\
&& \bfQ_\lambda[i,i] = \fr{2\lambda^2(t_i)}{h^2_{i-1}(h_{i-2} + h_{i-1})} + \fr{2\lambda^2(t_i)}{h_{i-1} h_i}\left(\fr{1}{h_{i-1}}+\fr{1}{h_i}\right) + \fr{2\lambda^2(t_i)}{h_i^2(h_i + h_{i+1})},
\ea
with corrected boundary entries
\ba
&& \bfQ_\lambda[1,1]=\fr{2\lambda^2(t_1)}{h_1^2(h_1+h_2)}, ~~\bfQ_\lambda[2,1] = -\fr{2\lambda(t_1)\lambda(t_2)}{h_1^2h_2}, \\
&&\bfQ_\lambda[2, 2]=\fr{2\lambda^2(t_2)}{h_2^2(h_2+h_3)} + \fr{2\lambda^2(t_2)}{h_1h_2}\left(\fr{1}{h_1} + \fr{1}{h_2}\right),\\
&& \bfQ_\lambda[n-1,n-1]=\fr{2\lambda^2(t_{n-1})}{h_{n-2}^2(h_{n-3}+h_{n-2})} + \fr{2\lambda^2(t_{n-1})}{h_{n-2}h_{n-1}}\left(\fr{1}{h_{n-2}} + \fr{1}{h_{n-1}}\right),\\
&&\bfQ_\lambda[n, n]=\fr{2\lambda^2(t_n)}{h_{n-1}^2(h_{n-2}+h_{n-1})}, ~~\bfQ_\lambda[n,n-1] = -\fr{2\lambda(t_{n-1})\lambda(t_n)}{h_{n-2}h_{n-1}^2}.
\ea

%

\subsection{Modeling adaptive smoothing function}

To implement fully Bayesian inference, we need a prior taken on the smoothing function $\lambda(t)$, which is assumed to be continuous and differentiable. Since it is restricted to be positive, we model $\lambda(t)$ on its log scale: $\nu(t) = \log(\lambda(t))$. \cite{yue:speck:non:10} and \cite{yue:spec:sun:aism:12} have proved that the prior on $\nu(t)$ must be proper in order to guarantee a proper posterior for such adaptive smoothing models. 

It is intuitive to model $\nu(t)$ in a similar way to $f(t)$. We therefore follow the basis expansion in (\ref{approx:f}) and represent $\nu(t)$ as a weighted sum of $m$ basis function $\omega_k(t)$, that is 
$$
\nu(t) = \sum_{k=1}^m\omega_k(t)\gamma_k,
$$
with random weights $\bfgamma=(\gamma_1,\ldots,\gamma_m)'$. Unfortunately, the previous GMRF prior cannot be put on $\bfgamma$ since it is intrinsic. \cite{lind:rue:11} derived an explicit link between GMRF and common Gaussian fields by considering SDE
\be\label{sde:lamb}
\left(\kappa^2 - d^2/dt^2\right)\nu(t) = dW(t)/dt,\quad
\ee
where $\kappa>0$ is fixed. Again, we use Galerkin method to weakly solve (\ref{sde:lamb}) as
\be\label{sde:weak:finite:lamb}
\sum_{k=1}^m \gamma_k\left\langle \omega_\ell, \kappa^2 - d^2\omega_k/dt^2\right\rangle \stackrel{d}{=} \left\langle \omega_\ell, dW/dt\right\rangle,\quad \ell=1,\ldots,m.
\ee
With piecewise linear basis, it can be shown that the left hand side of (\ref{sde:weak:finite:lamb}) is $(\kappa\bfB - \bfH)\bfgamma$ and the right hand side is a Gaussian random vector as before. As a result, the precision matrix of the corresponding GMRF is given by
\ba
\bfR =(\kappa^2\bfB - \bfH)'\bfB^{-1}(\kappa^2\bfB - \bfH) = \kappa^4\bfB - \kappa^2(\bfH'+\bfH) + \bfH'\bfB^{-1}\bfH.
\ea
To make $\bfR$ sparse, we replace $\bfB$ by $\tilde\bfB$ as before. This GMRF prior is proper and it is getting intrinsic as $\kappa$ goes to zero. Due to the computational advantage of GMRF, it is feasible to use full-rank basis expansion $(m=n)$ to make the method fully automatic.

\section{Posterior inference}\label{sec:infer}
The fully Bayesian inference requires the hyperpriors on parameters $\tau$, $\delta$ and $\eta$. We choose diffuse but proper gamma priors, i.e.\ $\mbox{Gamma}(\epsilon, \epsilon)$ for $\epsilon = 0.001$. Then, the joint posterior distribution of both adaptive smoothing spline models can be written as
$$
[\bfy\mid\bfw,\tau][\bfw\mid\delta,\bfgamma][\bfgamma\mid\eta][\tau][\delta][\eta].
$$
To obtain the posterior distribution, we here present two different approaches. They are simulation method via Monte Carlo Markov chain (MCMC) and approximation method based on integrated nested Laplace approximation (INLA). 

\subsection{MCMC approach}
Let $\bfPsi=\{\psi_j(t_i)\}_{i,j=1}^n$ and $\bfOmega=\{\omega_k(t_\ell)\}_{k,\ell=1}^m$ be the matrices of basis functions for $f(t)$ and $\nu(t)$, respectively. Then, the hierarchical models have the following structure: 
\ba
&&\bfy = \bfPsi\bfw + \bfvareps,\quad \bfvareps\sim N(\bfzero, \tau^{-1}\bfI),\\
&& [\bfw\mid\delta,\bflambda] \propto |\delta\bfQ_\lambda|^{1/2}_+\exp\left(-\fr{\delta}{2}\bfw'\bfQ_\lambda\bfw\right),\\
&& \log(\lambda_i) =  \nu_i,\quad \bfnu = \bfOmega\bfgamma,\\
&& [\bfgamma\mid\eta] \propto  |\eta\bfR|^{1/2}\exp\left(-\fr{\eta}{2}\bfgamma'\bfR\bfgamma\right),\\
&&\tau\sim\mbox{Gamma}(a_\tau,b_\tau), \\
&&\delta\sim\mbox{Gamma}(a_\delta,b_\delta),\\
&&\eta\sim\mbox{Gamma}(a_\eta,b_\eta).
\ea
We here focus on how to sample $\bfgamma$ from its full conditional because the rest sampling procedures are straightforward. As we can see, the full conditional of $\bfgamma$ is not a regular density, so we have to employ Metropolis-Hastings sampling technique. We here present an efficient algorithm to sample $\bfgamma$ when using the first adaptive SDE. Unfortunately, we have not found an equivalently efficient method for the second adaptive SDE, which, however, can be taken care of by INLA method as described in next section.   

A good proposal distribution is the key to the successful Metropolis-Hastings algorithm.
 It is helpful to see that the GMRF derived from the first adaptive SDE can be written as a random walk model, i.e.,
\ba
[\bfw\mid\delta,\bfgamma]\propto \prod_{i=1}^{n}\left(\delta e^{\gamma_i}\right)^{1/2}\exp\left(-\fr{\delta e^{\gamma_i}}{2}\tilde{w}_i^2\right),
\ea
where $\tilde\bfw = (0, \tilde w_2,\ldots,\tilde w_{n-1},0)' = \bfH\bfw$ (note the first and last rows of $\bfH$ are zeroes). Since $\gamma_i$ depends on $\tilde{w}_i$ only, it is possible to construct an accurate GMRF approximation for the full conditional of $\bfgamma$ given by $\mathcal{F}(\bfgamma\mid\bfw,\delta,\eta)\propto[\bfw\mid\delta,\bfgamma][\bfgamma\mid\eta]$ as follows. First, we approximate $[w_i\mid\delta,\gamma_i]$ using Taylor expansion at $\gamma_{0i}$, 
\ba
[w_i\mid\delta,\gamma_i]\approx\exp\left(a_i + b_i(\gamma_{0i})\gamma_i - \fr{1}{2}c_i(\gamma_{0i})\gamma_i^2\right),
\ea
where $a_i$ is the nuisance parameter, $b_i(\gamma_{0i}) = 1/2 - \delta e^{\gamma_{0i}}(1-\gamma_{0i})\tilde{w}_i^2/2$ and $c_i(\gamma_{0i}) = \delta e^{\gamma_{0i}}\tilde{w}_i^2/2 $. Letting $\bfb=(b_1,\ldots,b_n)'$ and $\bfc=(c_1,\ldots,c_n)'$, then the density
\ba
\mathcal{P}(\bfgamma\mid\bfgamma_0, \bfw,\delta,\eta)
&\propto& \exp\left(-\fr{\eta}{2}\bfgamma'\bfR\bfgamma +\sum_{i=1}^n\left(a_i + b_i(\gamma_{0i})\gamma_i - \fr{1}{2}c_i(\gamma_{0i})\gamma_i^2\right)\right)\\
&\propto&\exp\left(-\fr{1}{2}\bfgamma'(\eta\bfR + \mbox{diag}(\bfc))\bfgamma + \bfb'\bfgamma\right),
\ea 
is a GMRF approximation to $\mathcal{F}(\bfgamma\mid\bfw,\delta,\eta)$. In order to make the approximation accurate, we choose $\bfgamma_{0}$ to be the mode of $\mathcal{P}(\bfgamma\mid\bfgamma_0, \bfw,\delta,\eta)$, which can be obtained using, say Newton-Raphson method. Using the GMRF approximation as proposal distribution, we can update the whole $\bfgamma$ by accepting proposal $\bfgamma^*$ with probability 
\ba
\min \left (1, \quad\fr{\mathcal{F}(\bfgamma^*\mid\bfw,\delta,\eta)\mathcal{P}(\bfgamma\mid\bfgamma_0, \bfw,\delta,\eta)}{\mathcal{F}(\bfgamma\mid\bfw,\delta,\eta)\mathcal{P}(\bfgamma^*\mid\bfgamma_0, \bfw,\delta,\eta)}
\right).
\ea
Other full conditionals  are given by
\ba
&&(\bfw\mid\bfgamma,\delta,\tau)\sim N(\bfmu_w, \bfSigma_w), \quad\bfmu_w = \tau\bfSigma_w\bfPsi'\bfy \quad\mbox{and}\quad \bfSigma_w = \left(\tau\bfPsi'\bfPsi + \delta\bfQ_\lambda\right)^{-1}\\
&&(\tau\mid\bfw)\sim\mbox{Gamma}(n/2+a_\tau,\|\bfy-\bfPsi\bfw\|^2/2+b_\tau)\\
&& (\delta\mid\bflambda,\bfw)\sim\mbox{Gamma}(n/2-1+a_\delta,\bfw'\bfQ_\lambda\bfw/2+b_\delta)\\
&& (\eta\mid\bfgamma)\sim\mbox{Gamma}(n/2-1+a_\eta,\bfgamma'\bfR\bfgamma/2 + b_\eta),
\ea
all of which can be easily sampled.

\subsection{INLA approach}
\cite{Rue:Mart:Chop:inla:2009} have developed the \verb"R" computer package INLA for Bayesian inference using integrated nested Laplace approximations. The INLA can handle general Gaussian hierarchical models, including the both adaptive smoothing spline models developed in this paper. It accurately approximates marginal posterior densities and computes estimates much faster than general MCMC techniques. 

The general Gaussian hierarchical models have a set of hyperparameters $\bftheta$ with prior $\pi(\bftheta)$, a latent variable $\bff$ with density $\pi(\bff|\bftheta)$ and an observed response $\bfy$ with likelihood $\pi(\bfy|\bff,\bftheta)$. The posterior is then given by
\ba
\pi(\bff,\bftheta|\bfy)\propto\pi(\bfy|\bff,\bftheta)\pi(\bff|\bftheta)\pi(\bftheta).
\ea
We need to find the posterior marginals $\pi(f_i|\bfy)$ and $\pi(\theta_j|\bfy)$, which can be done using INLA. The approach is based on the following approximation for the posterior marginal of $\bftheta$:
\ba
\tilde{\pi}(\bftheta|\bfy) \propto \left.\fr{\pi(\bff,\bftheta,\bfy)}{\pi_G(\bff|\bftheta,\bfy)}\right|_{\bff=\bff^\star(\bftheta)},
\ea
where $\pi_G(\bff|\bftheta,\bfy)$ is the Gaussian approximation to the full conditional of $\bff$, and $\bff^\star(\bftheta)$ is the mode of the full conditional of $\bff$. 
The approximated marginals are then constructed as follows:
\ba
\tilde{\pi}(\theta_j|\bfy) &=& \int\tilde{\pi}(\bftheta|\bfy)d\bftheta_{-j},\\
\tilde{\pi}(f_i | \bfy) &=& \int \tilde{\pi}(f_i|\bftheta,\bfy)\tilde{\pi}(\bftheta|\bfy)d\bftheta,
\ea
where $\bftheta_{-j}$ denotes a subvector of $\bftheta$ without element $\theta_j$. The approximated marginal of $\theta_j$ can be obtained by summing out the remaining variables $\bftheta_{-j}$ from $\tilde{\pi}(\bftheta|\bfy)$. The approximated marginal of $f_i$ is obtained by, first, approximating the full conditional of $f_i$ with another Laplace approximation:
\ba
\tilde{\pi}(f_i | \bftheta,\bfy) \propto \left.\fr{\pi(\bff,\bftheta,\bfy)}{\pi_{GG}(\bff_{-i} | f_i,\bftheta,\bfy)}\right|_{\bff_{-i}=\bff^\star_{-i}(f_{-i},\bftheta)},
\ea
where $\tilde{\pi}_{GG}$ is the Gaussian approximation to $\bff_{-i} | f_i,\bftheta,\bfy$ and $\bff^\star_{-i}(f_{-i},\bftheta)$ is the mode configuration. Then, we numerically integrate out the parameters $\bftheta$ from $ \tilde{\pi}(f_i|\bftheta,\bfy)$. This nested approach makes the Laplace approximations very accurate. 

However, INLA has a limitation that is it only works when the number of hyperparameters in $\bftheta$ is small, say less than 15. The reason is that it becomes extremely expensive to numerically integrate out $\bftheta$ as its dimension increases. In our case, the hyperparameters $\bftheta = (\bfgamma, \tau, \delta,\eta)$. As a result, we have to use reduced-rank basis to model $\bfgamma$ if we want to fit the models with INLA.


\begin{figure}
\centering
\begin{tabular}{lll}
 \vspace{-.2in}
(a) Example 1 & (b) Example 2 & (c) Example 3\\  \vspace{-.2in}
\includegraphics[width=2.2in]{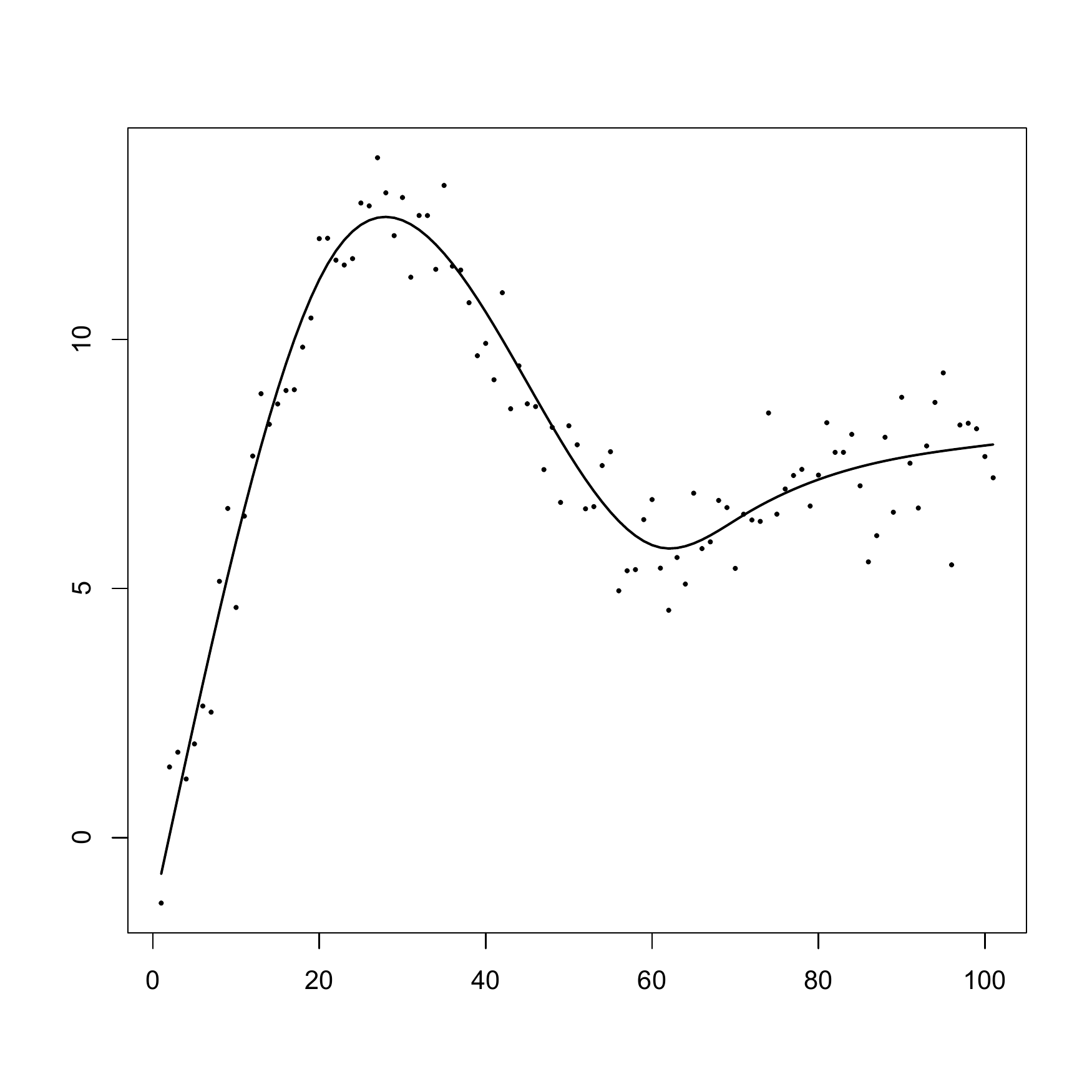} &
\includegraphics[width=2.2in]{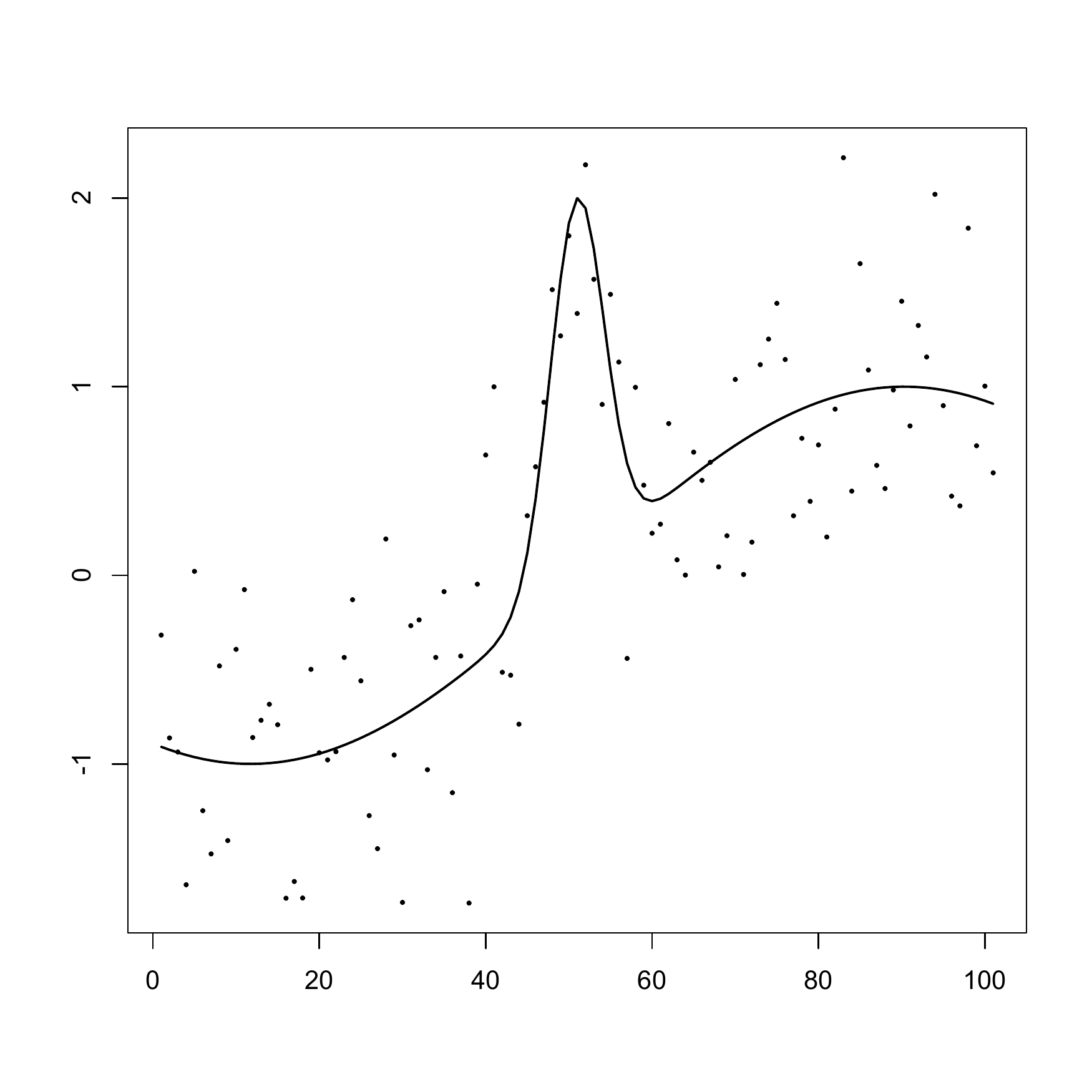} &
\includegraphics[width=2.2in]{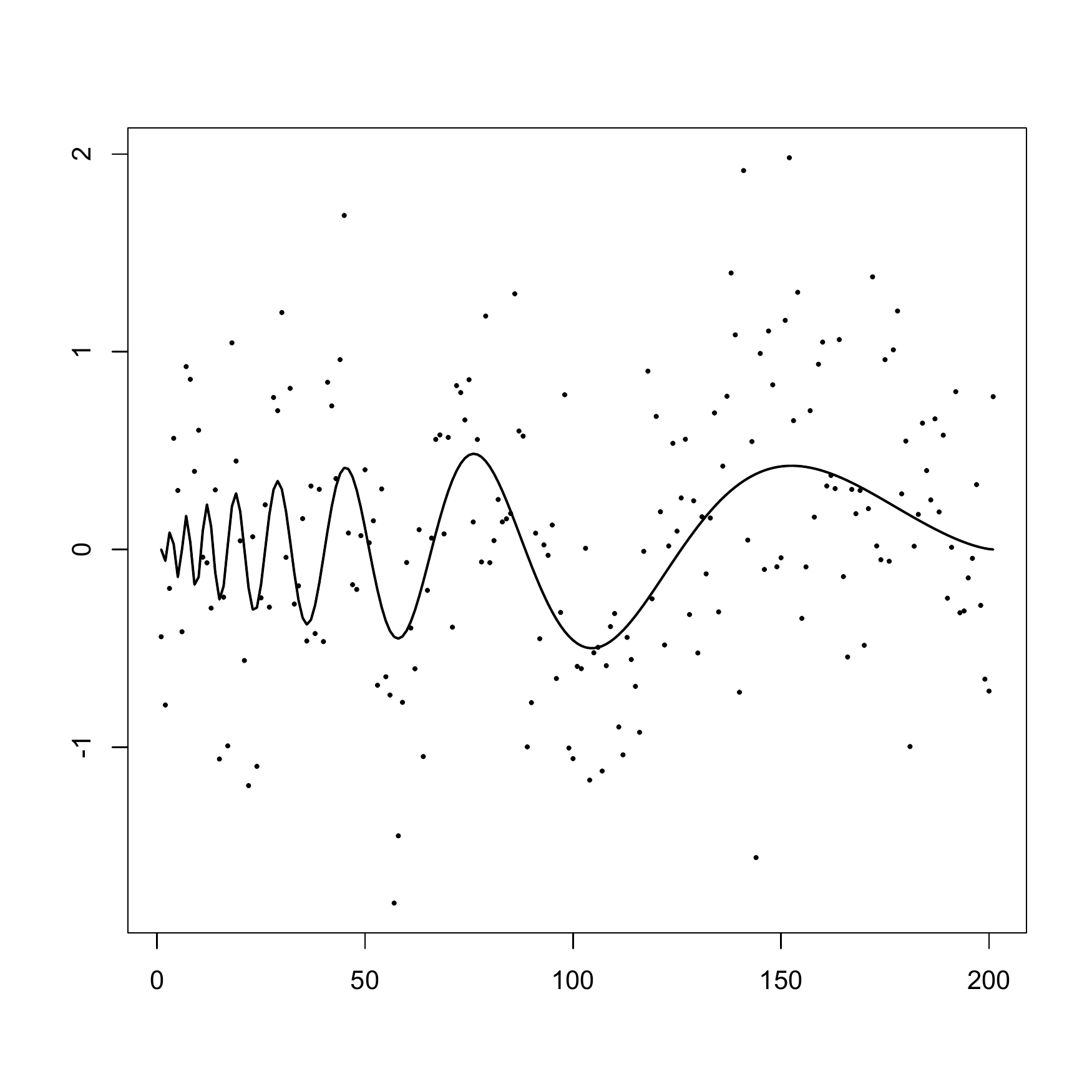} \\
\end{tabular}
\caption{The three true functions used in the simulation study together with one sample.}\label{fig:fun}
\end{figure}

\section{Simulated examples}\label{sec:simu}

In this section we consider three functions: a slowly-varying smooth function, a function with a sharp peak, that is spatially inhomogeneously smooth, and a highly-oscillating Doppler function. Gaussian noise is added to each in generating the data. The functions together with samples of data are shown in Figure \ref{fig:fun}. In Example 1, the true function is a spline with three internal knots at (0.2, 0.6, 0.7) and coefficients (20, 4, 6, 11, 6). The function is evaluated on a regular grid of 101 points, and a zero-mean Gaussian noise is added to the true function with standard deviation 0.9. In Example 2, the true function is $f(t) = \sin(t) + 2\exp(-30 t^2)$ for $t\in[-2,2]$, evaluated at 101 regularly spaced points, and the standard deviation of the noise is 0.5. In Example 3, the Doppler function is given by $f(t) = \sqrt{t(1-t)}\sin(2\pi(1+\epsilon)/(t+\epsilon))$ for $\epsilon=0.125$, evaluated at 201 regularly spaced points, and the standard deviation of the noise is 0.2.

We compare our Bayesian adaptive smoothing spline (BASS) estimates with ordinary smoothing spline (OSS) estimates, using mean squared error 
$$
\mbox{MSE} = \fr{1}{n}\sum_{i=1}^n\left[\hat f(t_i) - f(t_i)\right]^2.
$$
The BASS model derived from the first adaptive SDE is fitted by MCMC while the one from the second adaptive SDE is estimated by INLA.
Note that with MCMC we use the same number of knots as the data points, while with INLA we respectively use 3, 5 and 10 knots for the three examples. The median mean squared error, together with first and third quartile, based on 200 samples of data is reported in Table \ref{tab:mse}. As we can see, the OSS model slightly outperforms the two BASS models when estimating the slowly-varying smooth function, but the BASS models significantly work better in the peak and Doppler functions that are more spatially adaptive. It is interesting to see that two different BASS models, which are fitted by different methods, yield quite similar average MSE's. It indicates that the both SDE formulations offer appropriate adaptive smoothing, and INLA makes as accurate inference as MCMC does with much faster computation.   
The only limitation of INLA, as mentioned, is that it only works when there are a small number of hyperparameters to estimate. Therefore, INLA could be a better inferential tool than MCMC for BASS models if only a few knots are needed to capture the structure of the adaptive smoothing function. 

\begin{table}
\centering
\begin{tabular}{lccc}
\hline\hline
& Example 1 & Example 2 & Example 3\\
\hline
BASS-v1 &0.0620 (0.0444, 0.0854)&0.0297 (0.0219, 0.0420)&0.0072 (0.0061, 0.0084)\\
BASS-v2 &0.0633 (0.0468, 0.0886) &0.0274 (0.0206, 0.0377)&0.0072 (0.0058, 0.0088)\\
OSS &0.0600 (0.0401, 0.0853)&0.0408 (0.0336, 0.0531)&0.0092 (0.0082, 0.0102)\\
\hline
\end{tabular}
\caption{Simulation study. Median MSE with first and third quartiles in brackets based on 200 samples obtained using BASS and OSS procedures. The BASS-v1 and BASS-v2 denote the models derived from the first and second adaptive SDEs, respectively. }\label{tab:mse}
\end{table}


\section{Real data example}\label{sec:app}
\begin{figure}
\centering
\includegraphics[width=5in]{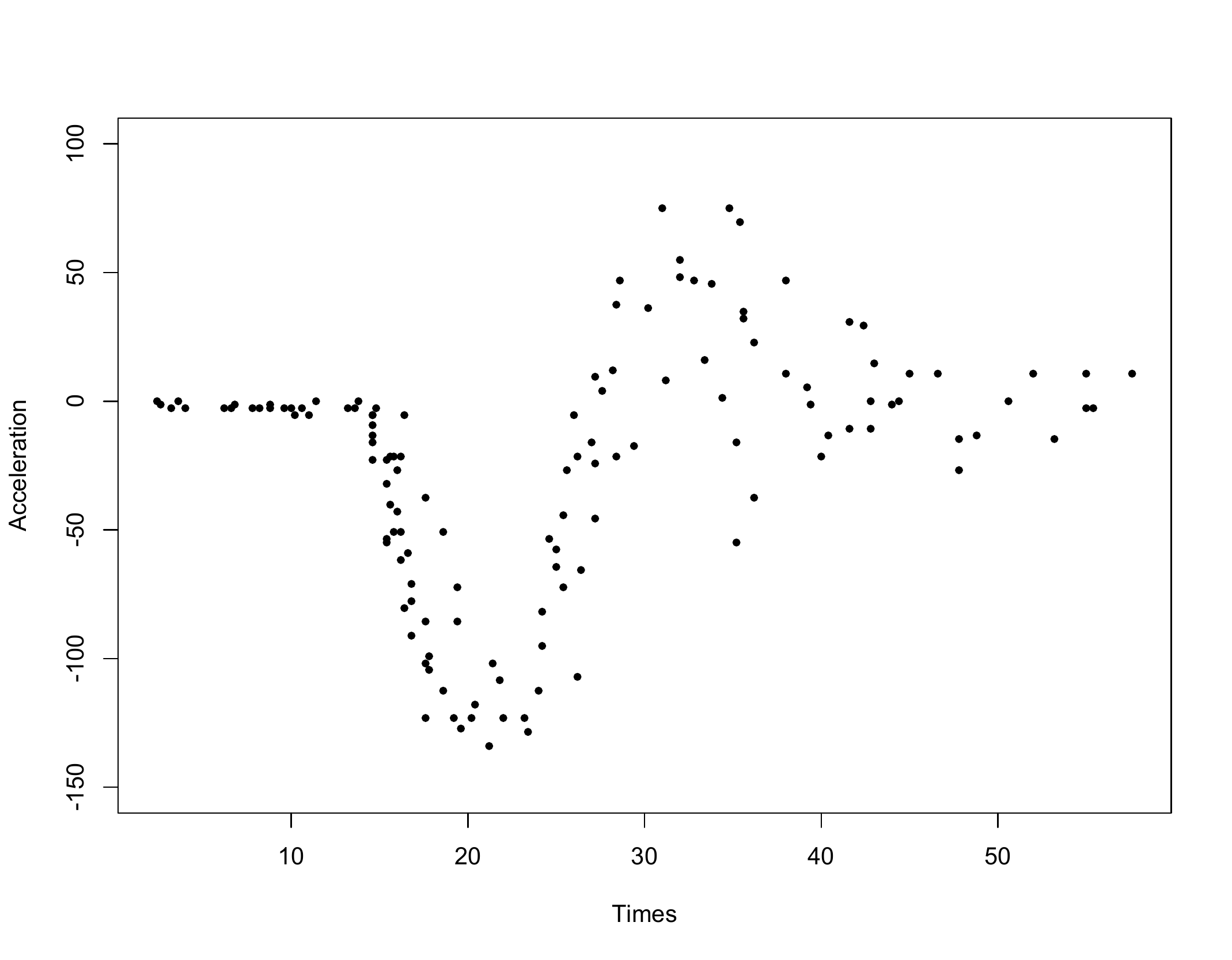} 
\vspace{-.3in}
\caption{The motorcycle impact data.}\label{fig:mcycle:data}
\end{figure}

To illustrate the techniques developed so far, we now consider the data presented in Figure \ref{fig:mcycle:data}. These observations consist of accelerometer readings taken through time in an experiment on the efficacy of crash helmets. 
The data set was used by \cite{silverman85} and is available in \verb@R@ software package. For various reasons, the time points are not regularly spaced, and there are multiple observations at some time points. In addition the observations are all subject to error. It is of interest both to discern the general shape of the underlying acceleration curve and to draw inferences about its minimum and maximum values. But, for illustrative purposes we shall concentrate on estimating the general shape only.

It is clear from Figure \ref{fig:mcycle:data} that the variance of the data is not constant over time. To take into account this heteroskedastic property, we modify model (\ref{mod:nonpar}) by adding random weights to the errors, that is $\vareps_i\sim N(0, \tau^{-1}\rho_i^{-1})$ for $i=1,\ldots,n$. Again, we take diffuse gamma prior on $\tau$. Regarding $\rho_i$, we use independent gamma prior with both shape and scale being half, i.e., $\rho_i \sim \mbox{Gamma}(0.5, 0.5)$. If integrate each $\rho_i$ out of $\vareps_i$, we can see that $\vareps_i$ follows an independent Cauchy distribution, which is able to provide flexible shrinkage due to its heavy tails and sharp peak. Such modifications on errors can be easily incorporated into the adaptive smoothing spline model by only adding the step of sampling $\rho_i$ to the MCMC algorithm.

\begin{figure}
\centering
\begin{tabular}{ll}
 \vspace{-.2in}
(a) OSS with Gaussian errors & (b) BASS with Gaussian errors \\ 
\includegraphics[width=3in]{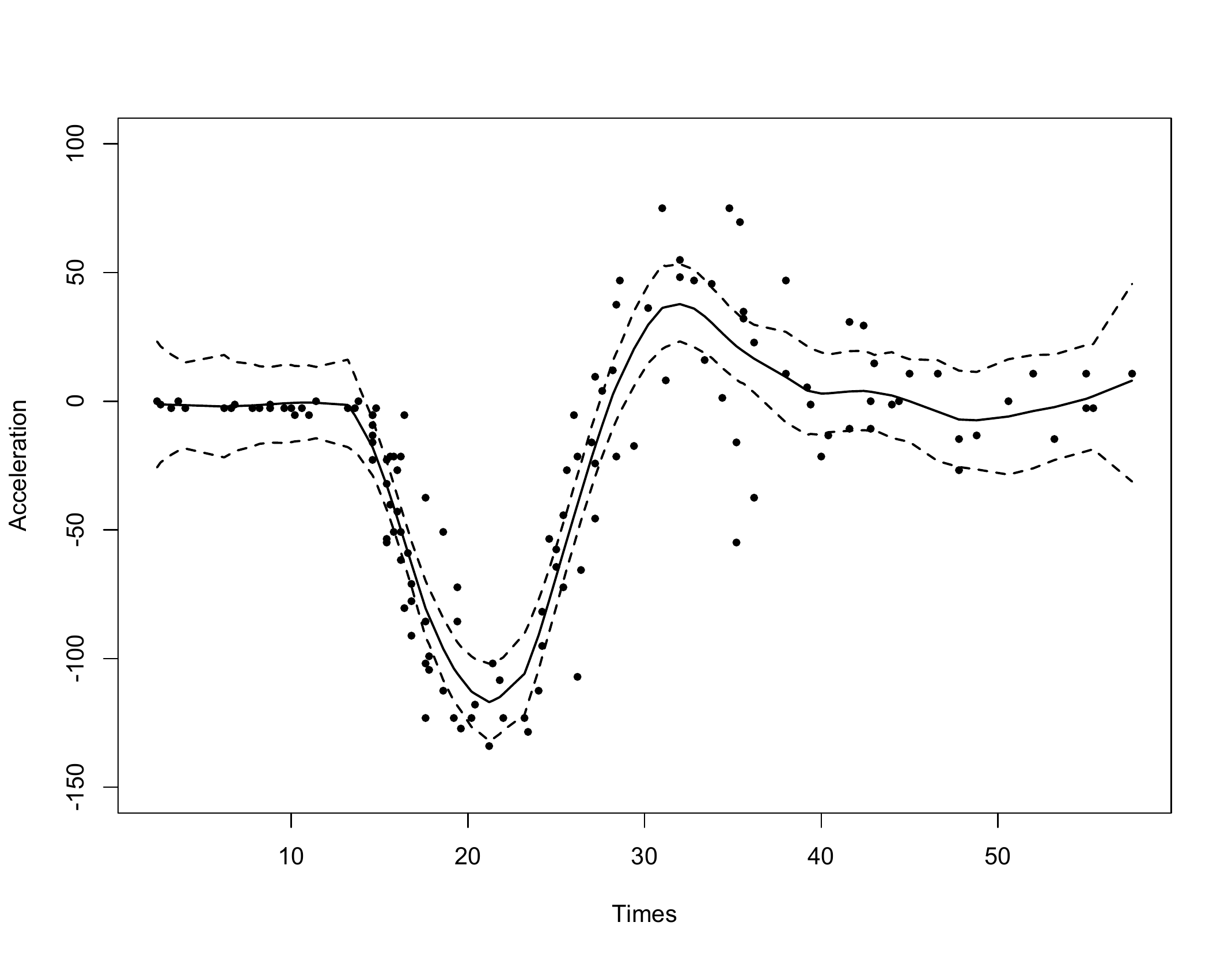} &
\includegraphics[width=3in]{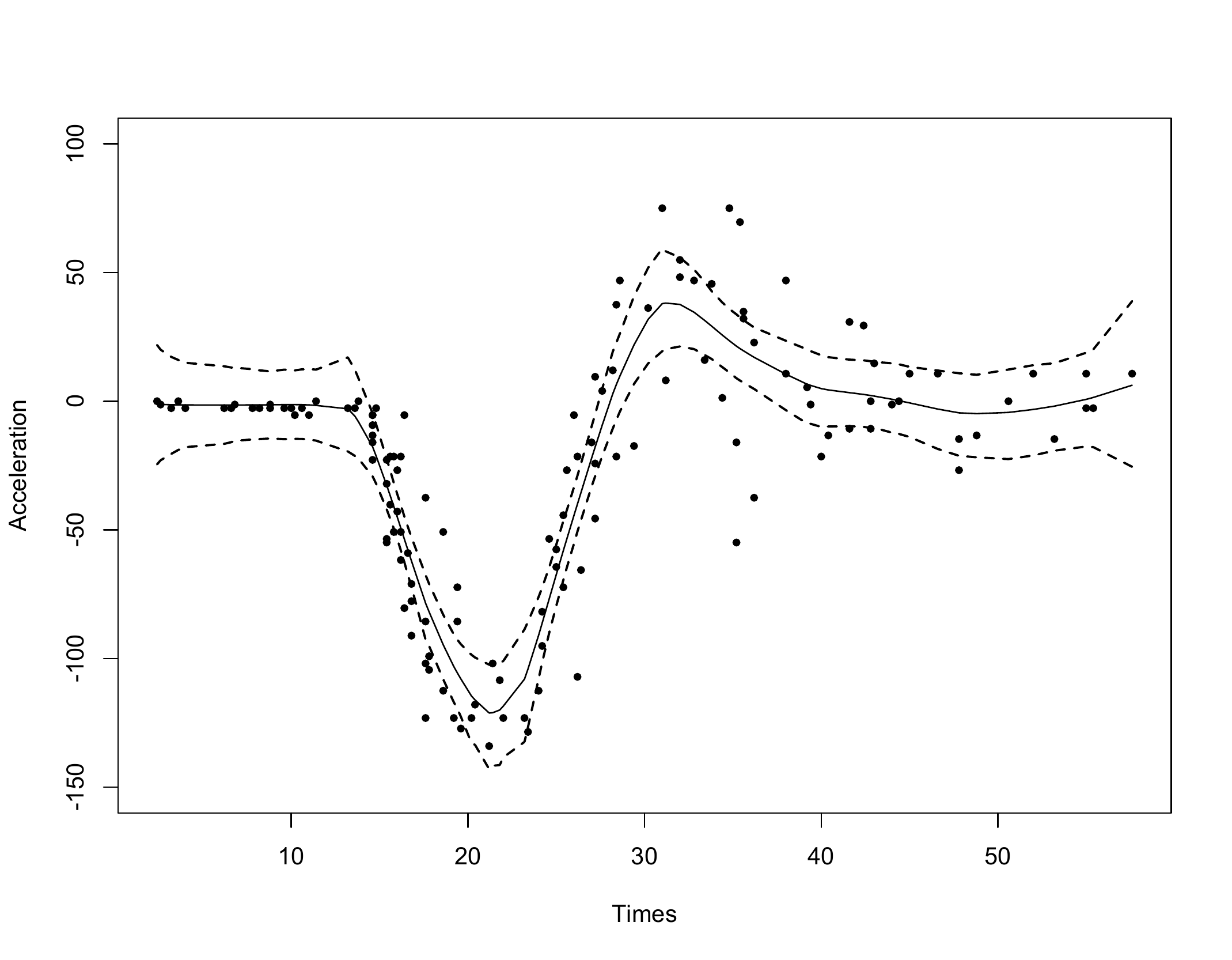} \\  \vspace{-.2in}
(c) OSS with Cauchy errors & (d) BASS with Cauchy errors \\
\includegraphics[width=3in]{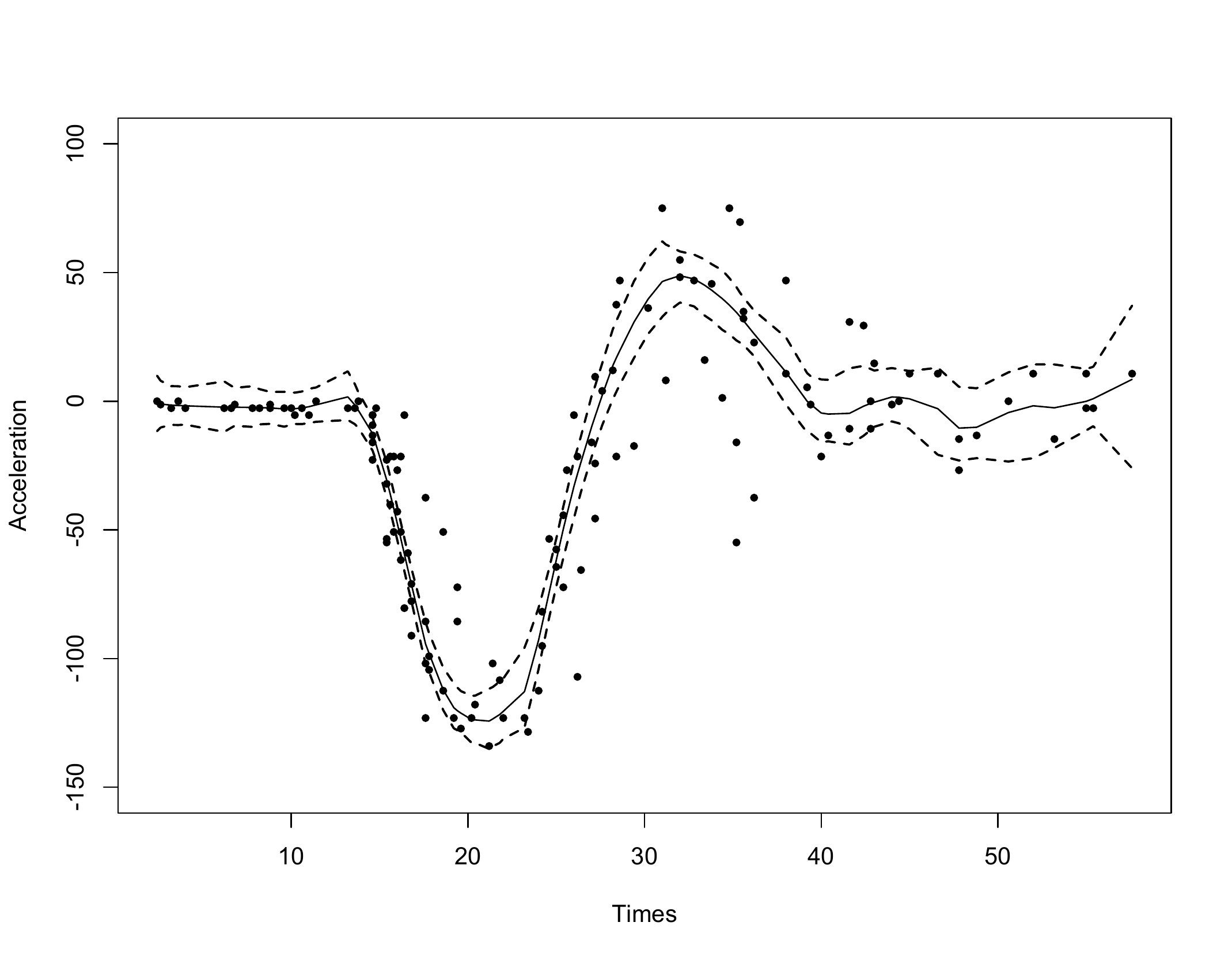} &
\includegraphics[width=3in]{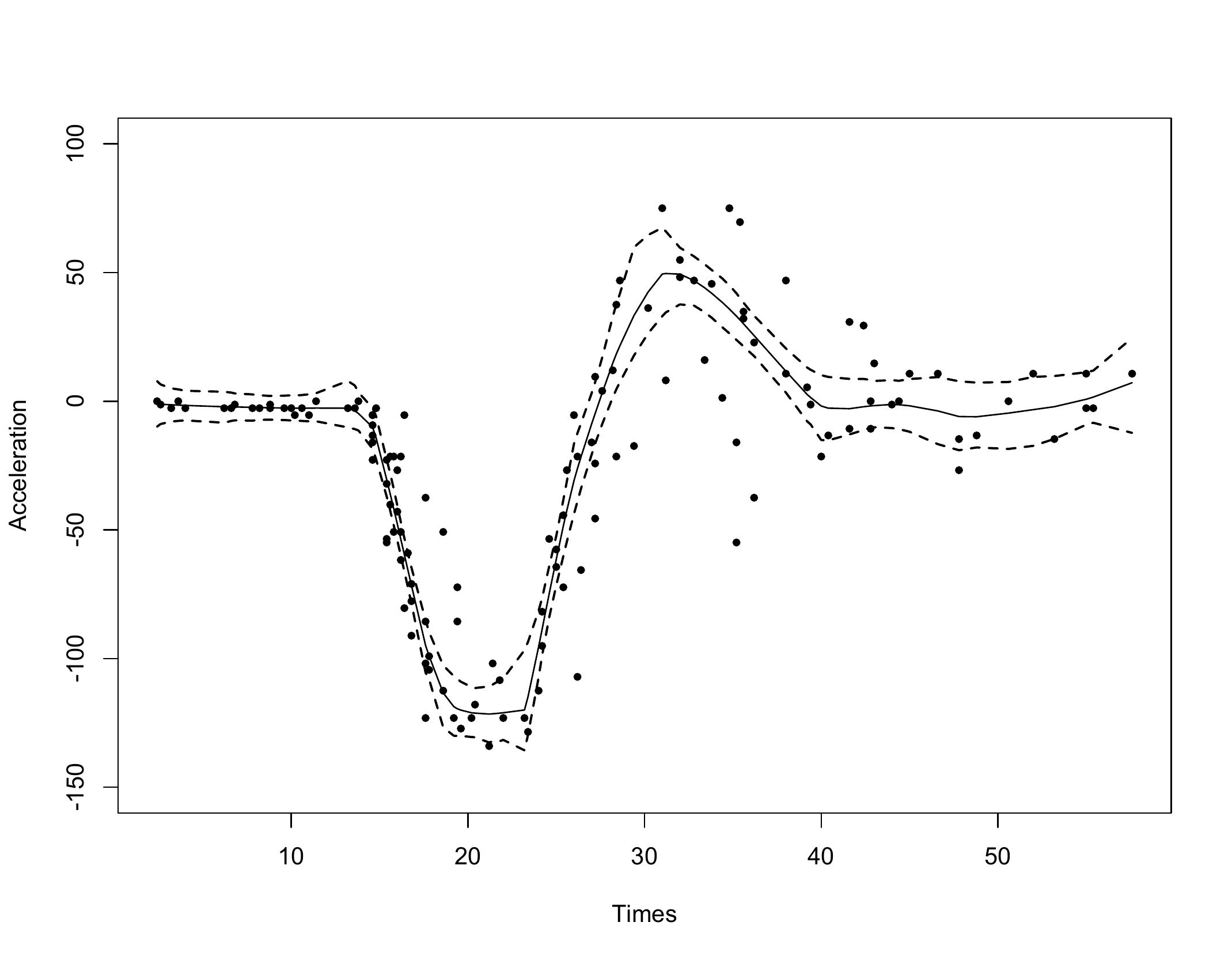} \\
\end{tabular}
\vspace{-.1in}
\caption{The fitted curved with their 95\% credible intervals, constructed from the motorcycle impact data.}\label{fig:mcycle:fit}
\end{figure}

The effect of applying adaptive smoothing technique and Cauchy errors is shown in Figure \ref{fig:mcycle:fit}. We here present four different fitted curves and their 95\% credible intervals: (a) OSS with Gaussian errors; (b) BASS with Gaussian errors; (c) OSS with Cauchy errors; (d) BASS with Cauchy errors. Note that we fit the BASS model derived from the first adaptive SDE using MCMC since the other model yields similar performance in the simulation study. 
As we can see, all the fits give a clear indication of the general pattern of the data, which is constant at first and then drops sharply, followed by a rebound above its original level before setting back. Compared to the OSS models, the BASS models show attractive adaptive smoothing features: the fits are smoother near the left and in the right half of the picture, while they yield lower drops in the middle. Compared to the Gaussian errors, the Cauchy errors make the fit follows the data more closely and offers a more reasonable credible interval (being narrow on the left end but wide in the right half), which captures the variance pattern well. In our opinion, the BASS model with Cauchy errors gives the best overall fit. 

%

\section{Conclusion}
In this paper we have developed a unified Bayesian approach to model adaptive smoothing splines. It is based on the connection between smoothing splines and stochastic differential equations. We showed that the SDE approach in \cite{lind:rue:08} can be easily adapted to adaptive smoothing problems. Using the finite element method, the GMRF representations of the adaptive smoothing splines were explicitly derived. Furthermore, we proposed efficient MCMC and INLA algorithms to make Bayesian inference. Finally, we demonstrated the effectiveness of our method through a simulation study and an application to the motorcycle data.

\begin{appendix}
\section{Appendix}
This section includes detailed proofs for the weak solutions of both adaptive SDEs. 

\subsection{Adaptive SDE I}
Using basis expansion (\ref{approx:f}), the adaptive SDE (\ref{sde:ad:weak}) becomes a linear equation system, whose left hand side can be written as $\bfH_\lambda\bfw$.  We here show how to derive the non-zero entries of matrix $\bfH_\lambda$. Using integration-by-parts, we have $[i,j]$th entry of $\bfH_\lambda$ as
\ba
\bfH_\lambda[i,j] &=& \Big\langle \psi_i(t), \lambda(t)\psi_j''(t)\Big\rangle\\
&=& \int_{t_1}^{t_n}\lambda(t) \psi_i(t) \psi_j''(t) dt \\
&=&\left.\lambda(t)\psi_i(t)\psi_j'(t)\right|_{t_1}^{t_n} - \int_{t_1}^{t_n}[\lambda(t)\psi_i(t)]'\psi_j'(t)dt\\
&=&\left.\lambda(t)\psi_i(t)\psi_j'(t)\right|_{t_1}^{t_n} - \int_{t_1}^{t_n}\lambda'(t)\psi_i(t)\psi_j'(t) dt - \int_{t_1}^{t_n} \lambda(t)\psi_i'(t)\psi_j'(t)dt.
\ea 
Since basis $\psi_i$ only overlap for neighboring locations, the nonzero entries in $i$th row of $\bfH_\lambda$ are $\bfH_\lambda[i,i-1]$, $\bfH_\lambda[i,i]$ and $\bfH_\lambda[i,i+1]$ for $i=2,\ldots,n-1$. Specifically, we have
\ba
\bfH[i,i-1] &=&  - \int_{t_{i-1}}^{t_i}\lambda'(t)\psi_i(t)\psi_{i-1}'(t) dt - \int_{t_{i-1}}^{t_i} \lambda(t)\psi_i'(t)\psi_{i-1}'(t)dt\\
&=& -\left.\lambda(t)\psi_i(t)\psi_{i-1}'(t)\right|_{t_{i-1}}^{t_i} + \int_{t_{i-1}}^{t_i}\lambda(t)[\psi_i(t)\psi_{i-1}'(t)]' dt - \int_{t_{i-1}}^{t_i} \lambda(t)\psi_i'(t)\psi_{i-1}'(t)dt\\
&=& -\left.\lambda(t)\psi_i(t)\psi_{i-1}'(t)\right|_{t_{i-1}}^{t_i} +  \int_{t_{i-1}}^{t_i} \lambda(t)\psi_i'(t)\psi_{i-1}'(t)dt - \int_{t_{i-1}}^{t_i} \lambda(t)\psi_i'(t)\psi_{i-1}'(t)dt\\
&=& -\lambda(t_i)\psi_i(t_i)\psi_{i-1}'(t_i) + \lambda(t_{i-1})\psi_i(t_{i-1})\psi_{i-1}'(t_{i-1})\\
&=& \lambda(t_i)/h_{i-1}.
\ea
Note that $\psi_{i-1}'(t)$ is constant between $t_{i-1}$ and $t_i$, and thus we have $[\psi_i(t)\psi_{i-1}'(t)]'  = \psi_i'(t)\psi_{i-1}'(t)$. Similarly, we have
\ba
\bfH_\lambda[i,i] &=&  - \int_{t_{i-1}}^{t_{i+1}}\lambda'(t)\psi_i(t)\psi_{i}'(t) dt - \int_{t_{i-1}}^{t_{i+1}} \lambda(t)\psi_i'(t)^2dt\\
&=&- \int_{t_{i-1}}^{t_i}\lambda'(t)\psi_i(t)\psi_{i}'(t) dt - \int_{t_{i-1}}^{t_i} \lambda(t)\psi_i'(t)^2dt \\
&& \quad\quad - \int_{t_i}^{t_{i+1}}\lambda'(t)\psi_i(t)\psi_{i}'(t) dt - \int_{t_i}^{t_{i+1}} \lambda(t)\psi_i'(t)^2dt\\
&=& -\left. \lambda(t)\psi_i(t)\psi_i'(t)\right|_{t_{i-1}}^{t_i} + \int_{t_{i-1}}^{t_i} \lambda(t)\psi_i'(t)^2dt - \int_{t_{i-1}}^{t_i} \lambda(t)\psi_i'(t)^2dt\\
&& \quad\quad -\left. \lambda(t)\psi_i(t)\psi_i'(t)\right|_{t_i}^{t_{i+1}}+ \int_{t_i}^{t_{i+1}} \lambda(t)\psi_i'(t)^2dt -  \int_{t_i}^{t_{i+1}} \lambda(t)\psi_i'(t)^2dt\\
&=& -\Big[\lambda(t_i)\psi_i(t_i) -  \lambda(t_{i-1})\psi_i(t_{i-1})\Big] /h_{i-1} +  \Big[ \lambda(t_{i+1})\psi_i(t_{i+1}) - \lambda(t_i)\psi_i(t_i) \Big] /h_{i}\\
&=& -\lambda(t_i)\left(\fr{1}{h_{i-1}} + \fr{1}{h_i}\right), \\
\bfH_\lambda[i, i+1] &=&  - \int_{t_i}^{t_{i+1}}\lambda'(t)\psi_i(t)\psi_{i+1}'(x) dt - \int_{t_i}^{t_{i+1}} \lambda(t)\psi_i'(t)\psi_{i+1}'(t)dt\\
&=&-\left.\lambda(t)\psi_i(t)\psi_{i+1}'(t)\right|_{t_i}^{t_{i+1}} +  \int_{t_i}^{t_{i+1}} \lambda(t)\psi_i'(t)\psi_{i+1}'(t)dt - \int_{t_i}^{t_{i+1}} \lambda(t)\psi_i'(t)\psi_{i+1}'(t)dt\\
&=&-\lambda(t_{i+1})\psi_i(t_{i+1})\psi_{i+1}'(t_{i+1}) + \lambda(t_i)\psi_i(t_i)\psi_{i+1}'(t_i)\\
&=&  \lambda(t_i)/h_{i}.
\ea
For first and last row of $\bfH_\lambda$, the (possible) nonzero entries are $\bfH_\lambda[1,1]$, $\bfH_\lambda[1,2]$, $\bfH_\lambda[n,n-1]$ and $\bfH_\lambda[n,n]$, which happen to be zeroes due to the intrinsic condition. We here only show the derivation of $\bfH_\lambda[1,1]$ and $\bfH_\lambda[1,2]$, and the other entries can be obtained similarly. We have
\ba
\bfH_\lambda[1,1] &=& \left.\lambda(t)\psi_1(t)\psi_1'(t)\right|_{t_1}^{t_n} - \int_{t_1}^{t_2}\lambda'(t)\psi_1(t)\psi_1'(t) dt - \int_{t_1}^{t_2} \lambda(t)\psi_1'(t)^2dt\\
&=&-\lambda(t_1)\psi_1(t_1)\psi_1'(t_1) -  \left.\lambda(t)\psi_1(t)\psi_1'(t)\right|_{t_1}^{t_2} \\
&=&-\lambda(t_1)\psi_1(t_1)\psi_1'(t_1) +\lambda(t_1)\psi_1(t_1)\psi_1'(t_1) \\
&=& 0,\\
\bfH_\lambda[1,2] &=& \left.\lambda(t)\psi_1(t)\psi_2'(t)\right|_{t_1}^{t_n} - \int_{t_1}^{t_2}\lambda'(t)\psi_1(t)\psi_2'(t) dt - \int_{t_1}^{t_2} \lambda(t)\psi_1'(t)\psi_2'(t)dt\\
&=&-\lambda(t_1)\psi_1(t_1)\psi_2'(t_1) -  \left.\lambda(t)\psi_1(t)\psi_2'(t)\right|_{t_1}^{t_2} \\
&=&-\lambda(t_1)\psi_1(t_1)\psi_2'(t_1) +\lambda(t_1)\psi_1(t_1)\psi_2'(t_1) \\
&=& 0.
\ea
Finally, we can easily see that $\bfH_\lambda = \bfLambda\bfH$, where $\bfLambda$ is the diagonal matrix of $\lambda(\cdot)$'s and $\bfH$ is the tridiagonal matrix defined as in (\ref{sde:lhs}).

\subsection{Adaptive SDE II}
Letting $\tilde f = \lambda(t)f(t)$, the left hand side of (\ref{sde:weak2}) can be written as 
\ba
\Big\langle \psi_i(t), \tilde f''(t)\Big\rangle &=& \int \psi_i(t)\tilde f''(t) dt \\
&=&\left. \psi_i(t)\tilde f'(t)\right |_{t_1}^{t_n} - \int\psi_i'(t)\tilde f'(t) dt\\
&=&\left.\psi_i(t)\Big[\lambda'(t)f(t) + \lambda(t)f'(t)\Big]\right |_{t_1}^{t_n}- \int\psi_i'(t)\Big[\lambda'(t)f(t) + \lambda(t)f'(t)\Big] dt
\ea
Using basis expansion (\ref{approx:f}), the adaptive SDE (\ref{sde:weak2}) becomes a linear equation system, whose left hand side can be written as $\bfH_\lambda\bfw$. We then have $[i,j]$th entry of $\bfH_\lambda$ as
\ba
\bfH_\lambda[i,j] =\left. \lambda'(t)\psi_i(t)\psi_j(t)\right |_{t_1}^{t_n} + \left. \lambda(t)\psi_i(t)\psi'_j(t)\right |_{t_1}^{t_n} - \int_{t_1}^{t_n}\lambda'(t)\psi_i'(t)\psi_j(t)dt- \int_{t_1}^{t_n}\lambda(t)\psi_i'(t)\psi_j'(t)dt.
\ea 
Since basis $\psi_i$ only overlap for neighboring locations, the nonzero entries in $i$th row of $\bfH_\lambda$ are $\bfH_\lambda[i,i-1]$, $\bfH_\lambda[i,i]$ and $\bfH_\lambda[i,i+1]$ for $i=2,\ldots,n-1$. Specifically, we have
\ba
\bfH[i,i-1] &=&  - \int_{t_{i-1}}^{t_i}\lambda'(t)\psi_i'(t)\psi_{i-1}(t) dt - \int_{t_{i-1}}^{t_i} \lambda(t)\psi_i'(t)\psi_{i-1}'(t)dt\\
&=& -\left.\lambda(t)\psi_i'(t)\psi_{i-1}(t)\right|_{t_{i-1}}^{t_i} + \int_{t_{i-1}}^{t_i}\lambda(t)[\psi_i'(t)\psi_{i-1}(t)]' dt - \int_{t_{i-1}}^{t_i} \lambda(t)\psi_i'(t)\psi_{i-1}'(t)dt\\
&=& -\left.\lambda(t)\psi_i'(t)\psi_{i-1}(t)\right|_{t_{i-1}}^{t_i} +  \int_{t_{i-1}}^{t_i} \lambda(t)\psi_i'(t)\psi_{i-1}'(t)dt - \int_{t_{i-1}}^{t_i} \lambda(t)\psi_i'(t)\psi_{i-1}'(t)dt\\
&=& -\lambda(t_i)\psi_i'(t_i)\psi_{i-1}(t_i) + \lambda(t_{i-1})\psi_i'(t_{i-1})\psi_{i-1}(t_{i-1})\\
&=& \lambda(t_{i-1})/h_{i-1}.
\ea
Note that $\psi_{i}'(t)$ is constant between $t_{i-1}$ and $t_i$, and thus we have $[\psi_i'(t)\psi_{i-1}(t)]'  = \psi_i'(t)\psi_{i-1}'(t)$. Similarly, we have
\ba
\bfH_\lambda[i,i] &=&  - \int_{t_{i-1}}^{t_{i+1}}\lambda'(t)\psi_i'(t)\psi_{i}(t) dt - \int_{t_{i-1}}^{t_{i+1}} \lambda(t)\psi_i'(t)^2dt= -\lambda(t_i)\left(\fr{1}{h_{i-1}} + \fr{1}{h_i}\right), 
\ea
which is the same as in the previous case, and
\ba
\bfH_\lambda[i, i+1] &=&  - \int_{t_i}^{t_{i+1}}\lambda'(t)\psi_i'(t)\psi_{i+1}(t) dt - \int_{t_i}^{t_{i+1}} \lambda(t)\psi_i'(t)\psi_{i+1}'(t)dt\\
&=&-\left.\lambda(t)\psi_i'(t)\psi_{i+1}(t)\right|_{t_i}^{t_{i+1}} +  \int_{t_i}^{t_{i+1}} \lambda(t)\psi_i'(t)\psi_{i+1}'(t)dt - \int_{t_i}^{t_{i+1}} \lambda(t)\psi_i'(t)\psi_{i+1}'(t)dt\\
&=&-\lambda(t_{i+1})\psi_i'(t_{i+1})\psi_{i+1}(t_{i+1}) + \lambda(t_i)\psi_i'(t_i)\psi_{i+1}(t_i)\\
&=&  \lambda(t_{i+1})/h_{i}.
\ea
For first and last row of $\bfH_\lambda$, the (possible) nonzero entries are $\bfH_\lambda[1,1]$, $\bfH_\lambda[1,2]$, $\bfH_\lambda[n,n-1]$ and $\bfH_\lambda[n,n]$, of which the first two entries can be derived as
\ba
\bfH_\lambda[1,1] &=&\left.\lambda'(t)\psi_1^2(t)\right|_{t_1}^{t_n} +
 \left.\lambda(t)\psi_1(t)\psi_1'(t)\right|_{t_1}^{t_n} - \int_{t_1}^{t_2}\lambda'(t)\psi'_1(t)\psi_1(t) dt - \int_{t_1}^{t_2} \lambda(t)\psi_1'(t)^2dt\\
&=&-\lambda'(t_1)\psi_1^2(t_1) - \lambda(t_1)\psi_1'(t_1)\psi_1(t_1) - \left.\lambda(t_1)\psi_1'(t)\psi_1(t)\right |_{t_1}^{t_2}\\
&=& -\lambda'(t_1)\psi_1^2(t_1) - \lambda(t_1)\psi_1'(t_1)\psi_1(t_1) + \lambda(t_1)\psi_1'(t_1)\psi_1(t_1)\\
&=& -\lambda'(t_1),\\
\bfH_\lambda[1,2] &=& \left.\lambda'(t)\psi_1(t)\psi_2(t)\right|_{t_1}^{t_n} +
 \left.\lambda(t)\psi_1(t)\psi_2'(t)\right|_{t_1}^{t_n} - \int_{t_1}^{t_2}\lambda'(t)\psi'_1(t)\psi_2(t) dt - \int_{t_1}^{t_2} \lambda(t)\psi_1'(t)\psi_2'(t)dt\\
&=&- \lambda'(t_1)\psi_1(t_1)\psi_2(t_1) - \lambda(t_1)\psi_1(t_1)\psi_2'(t_1) -\lambda(t_{2})\psi_1'(t_{2})\psi_{2}(t_{2}) + \lambda(t_1)\psi_1'(t_1)\psi_{2}(t_1)\\
&=& -\lambda(t_1)/h_1 + \lambda(t_2)/h_1.
\ea
Similarly, the last two entries are given by 
$$
\bfH_\lambda[n-1,n] = \lambda(t_{n-1})/h_{n-1} - \lambda(t_n)/h_{n-1}\quad\mbox{and}\quad \bfH_{\lambda}[n,n] = \lambda'(t_n).
$$
These four entries can be viewed as (at least approximately) the derivatives of $\lambda(t)$ at the boundary points. To be consistent with the previous case, we assume the Neumann boundary condition: $\lambda'(t_1) = \lambda'(t_n) = 0$, to make the entries be zeroes. Then, we can easily see that $\bfH_\lambda = \bfH\bfLambda$, where $\bfLambda$ is the diagonal matrix of $\lambda(\cdot)$'s and $\bfH$ is the matrix defined as in (\ref{sde:lhs}). 

\end{appendix}

\bigskip
\renewcommand{\baselinestretch}{1.0}
\tiny\normalsize
\bibliography{yuyuestats}
\bibliographystyle{sinica}

\end{document}

%% file: Macros.tex
\newcommand{\nmathbf}{\bm}

\def\bfA{\nmathbf A}
\def\bfB{\nmathbf B}

\def\bfH{\nmathbf H}
\def\bfI{\nmathbf I}

\def\bfQ{\nmathbf Q}
\def\bfR{\nmathbf R}

\def\bfb{\nmathbf b}
\def\bfc{\nmathbf c}

\def\bff{\nmathbf f}

\def\bfw{\nmathbf w}

\def\bfy{\nmathbf y}

\def\bfgamma  {\nmathbf \gamma}

\def\bftheta  {\nmathbf \theta}

\def\bflambda {\nmathbf \lambda}
\def\bfmu     {\nmathbf \mu}
\def\bfnu     {\nmathbf \nu}

\def\bfLambda {\nmathbf \Lambda}

\def\bfPsi     {\nmathbf \Psi}

\def\bfSigma  {\nmathbf \Sigma}

\def\bfOmega  {\nmathbf \Omega}

\newcommand{\bfzero}{{\nmathbf 0}}

\newcommand{\vareps}{\varepsilon}
\def\bfvareps{\nmathbf \varepsilon}

\def\boldfacefake#1{\kern-4pt
    \hbox{ \mathsurround=0pt
    \hbox to 0.4pt{$#1$\hss}\hbox to 0.4pt{$#1$\hss}\hbox {$#1$}}}

\newcommand{\bpro}{\begin{proof}}
\newcommand{\epro}{\end{proof}}
\newcommand{\be}{\begin{eqnarray}}
\newcommand{\ee}{\end{eqnarray}}
\newcommand{\ba}{\begin{eqnarray*}}
\newcommand{\ea}{\end{eqnarray*}}
\newcommand{\bc}{\begin{center}}
\newcommand{\ec}{\end{center}}
\newcommand{\btab}[1]{\begin{tabular}{#1}}
\newcommand{\etab}{\end{tabular}}

\newcommand{\fr}{\frac}

\newtheorem{theorem0}{Theorem}
\newtheorem{lemma0}{Lemma}
\newtheorem{remark0}{Remark}
\newtheorem{fact0}{Fact}
\newtheorem{example0}{Example}
\newtheorem{definition0}{Definition}
\newtheorem{corollary0}{Corollary}
\newtheorem{proposition0}{Proposition}
\newtheorem{algorithmY}{Algorithm}

\newcommand{\reals}{\mbox{\rm I\kern-.20em R}}
\newcommand{\sreals}{\mbox{\small \rm I\kern-.20em R}}

\newcommand{\goto}{\rightarrow}

\newenvironment{proof}{\vspace*{-2mm}\noindent{\sc Proof}. \mbox{} \rm}
   {\hfill\fbox{}}

%% file: ss_spde_v3.bbl
\begin{thebibliography}{44}
\expandafter\ifx\csname natexlab\endcsname\relax\def\natexlab#1{#1}\fi

\bibitem[{Abramovich and Steinberg(1996)}]{Abra:Stei:impr:1996}
Abramovich, F. and Steinberg, D.~M. (1996). Improved inference in nonparametric
  regression using ${L}_k$-smoothing splines. {\em Journal of Statistical
  Planning and Inference\/} {\bf 49}, 327--341.

\bibitem[{Baladandayuthapani et~al.(2005)Baladandayuthapani, Mallick and
  Carroll}]{Bala:Mall:Carr:spat:2005}
Baladandayuthapani, V., Mallick, B.~K. and Carroll, R.~J. (2005). Spatially
  adaptive {B}ayesian penalized regression splines ({P}-splines). {\em Journal
  of Computational and Graphical Statistics\/} {\bf 14}, 378--394.

\bibitem[{Brezger and Lang(2006)}]{brezger:csda:06}
Brezger, A. and Lang, S. (2006). Generalized structured additive regression
  based on {B}ayesian {P}-splines. {\em Computational Statistics and Data
  Analysis\/} {\bf 50}, 967--991.

\bibitem[{Crainiceanu et~al.(2007)Crainiceanu, Ruppert, Carroll, Adarsh and
  Goodner}]{Crai:Rupp:07}
Crainiceanu, C., Ruppert, D., Carroll, R., Adarsh, J. and Goodner, B. (2007).
  Spatially adaptive {P}enalized splines with heteroscedastic errors. {\em
  Journal of Computational and Graphical Statistics\/} 265--288.

\bibitem[{Cummins et~al.(2001)Cummins, Filloon and
  Nychka}]{Cumm:Fill:Nych:2001}
Cummins, D.~J., Filloon, T.~G. and Nychka, D. (2001). Confidence intervals for
  nonparametric curve estimates: {T}oward more uniform pointwise coverage. {\em
  Journal of the American Statistical Association\/} {\bf 96}, 233--246.

\bibitem[{Denison et~al.(1998)Denison, Mallick and Smith}]{denison98}
Denison, D. G.~T., Mallick, B.~K. and Smith, A. F.~M. (1998). Automatic
  bayesian curve fitting. {\em Journal of the Royal Statistical Society: Series
  B (Statistical Methodology)\/} {\bf 60}, 333--350.

\bibitem[{Di~Matteo et~al.(2001)Di~Matteo, Genovese and
  Kass}]{Dim:Gen:Kas:baye:2001}
Di~Matteo, I., Genovese, C.~R. and Kass, R.~E. (2001). Bayesian curve-fitting
  with free-knot splines. {\em Biometrika\/} {\bf 88}, 1055--1071.

\bibitem[{Eilers and Marx(1996)}]{Eile:Marx:pspline:1996}
Eilers, P. and Marx, B. (1996). Flexible smoothing with {B}-splines and
  penalties (with discussion). {\em Statistical Science\/} {\bf 11}, 89--121.

\bibitem[{Eilers and Marx(2010)}]{Eilers10}
Eilers, P. H.~C. and Marx, B.~D. (2010). Splines, knots, and penalties. {\em
  Wiley Interdisciplinary Reviews: Computational Statistics\/} {\bf 2},
  637--653.

\bibitem[{Eubank(1999)}]{Euba:nonp:1999}
Eubank, R.~L. (1999). {\em Nonparametric Regression and Spline Smoothing\/}.
  Marcel Dekker Inc.

\bibitem[{Fahrmeir and Knorr-Held(2000)}]{fahr:knor:00}
Fahrmeir, L. and Knorr-Held, L. (2000). Dynamic and semiparametric models. In
  {\em Smoothing and regression: approaches, computation, and application\/}
  (M.~G. Schimek, ed.), 513--544, New York: Wiley.

\bibitem[{Fahrmeir and Lang(2001)}]{Fahr:Lang:baye:2001}
Fahrmeir, L. and Lang, S. (2001). Bayesian inference for generalized additive
  mixed models based on {M}arkov random field priors. {\em Journal of the Royal
  Statistical Society, Series C: Applied Statistics\/} {\bf 50}, 201--220.

\bibitem[{Fahrmeir and Wagenpfeil(1996)}]{Fahr:Wage:smoo:1996}
Fahrmeir, L. and Wagenpfeil, S. (1996). Smoothing hazard functions and
  time-varying effects in discrete duration and competing risks models. {\em
  Journal of the American Statistical Association\/} {\bf 91}, 1584--1594.

\bibitem[{Fan and Gijbels(1996)}]{fan96}
Fan, J. and Gijbels, I. (1996). {\em Local Polynomial Modeling and its
  Application\/}. London: Chapman and Hall.

\bibitem[{Green and Silverman(1994)}]{Gree:Silv:nonp:1994}
Green, P.~J. and Silverman, B.~W. (1994). {\em Nonparametric Regression and
  Generalized Linear Models: a Roughness Penalty Approach\/}. Chapman \& Hall
  Ltd.

\bibitem[{Gu(2002)}]{Gu:smoo:2002}
Gu, C. (2002). {\em Smoothing Spline {ANOVA} Models\/}. Springer-Verlag Inc,
  New York.

\bibitem[{Hansen and Kooperberg(2002)}]{hansen02}
Hansen, M.~H. and Kooperberg, C. (2002). Spline adaptation in extended linear
  models (with discussion). {\em Statistical Science\/} {\bf 17}, 2--51.

\bibitem[{Holmes and Mallick(2001)}]{Holm:Mall:bay:2001}
Holmes, C.~C. and Mallick, B.~K. (2001). Bayesian regression with multivariate
  linear splines. {\em Journal of the Royal Statistical Society, Series B:
  Statistical Methodology\/} {\bf 63}, 3--17.

\bibitem[{Kimeldorf and Wahba(1970)}]{kim:wah:70}
Kimeldorf, G.~S. and Wahba, G. (1970). A correspondence between {B}ayesian
  estimation on stochastic processes and smoothing by splines. {\em Annals of
  Mathematical Statistics\/} {\bf 41}, 495--502.

\bibitem[{Krivobokova et~al.(2008)Krivobokova, Crainiceanu and
  Kauermann}]{Kriv:Crai:Kaue:fast:2008}
Krivobokova, T., Crainiceanu, C.~M. and Kauermann, G. (2008). Fast {A}daptive
  {P}enalized {S}plines. {\em Journal of Computational and Graphical
  Statistics\/} {\bf 17}, 1--20.

\bibitem[{Lang and Brezger(2004)}]{Lang:Brez:baye:2004}
Lang, S. and Brezger, A. (2004). Bayesian {P}-splines. {\em Journal of
  Computational and Graphical Statistics\/} {\bf 13}, 183--212.

\bibitem[{Lang et~al.(2002)Lang, Fronk and Fahrmeir}]{Lang:Fron:Fahr:func:2002}
Lang, S., Fronk, E.~M. and Fahrmeir, L. (2002). Function estimation with
  locally adaptive dynamic models. {\em Computational Statistics\/} {\bf 17},
  479--499.

\bibitem[{Lindgren and Rue(2008)}]{lind:rue:08}
Lindgren, F. and Rue, H. (2008). On the second-order random walk model for
  irregular locations. {\em Scandinavian Journal of Statistics\/} {\bf 35},
  691--700.

\bibitem[{Lindgren et~al.(2011)Lindgren, Rue and Lindstr{\"o}m}]{lind:rue:11}
Lindgren, F., Rue, H. and Lindstr{\"o}m, J. (2011). An explicit link between
  gaussian fields and gaussian markov random fields: the stochastic partial
  differential equation approach (with discussion). {\em Journal of the Royal
  Statistical Society: Series B (Statistical Methodology)\/} {\bf 73},
  423--498.

\bibitem[{Luo and Wahba(1997)}]{luo:wahba:97}
Luo, Z. and Wahba, G. (1997). Hybrid adaptive splines. {\em Journal of the
  American Statistical Association\/} {\bf 92}, 107--116.

\bibitem[{O'Sullivan(1986)}]{OSullivan86}
O'Sullivan, F. (1986). A statistical perspective on ill-posed inverse problems.
  {\em Statistical Science\/} {\bf 1}, 502--527.

\bibitem[{Pintore et~al.(2006)Pintore, Speckman and
  Holmes}]{Pin:Spec:Holm:spat:2006}
Pintore, A., Speckman, P.~L. and Holmes, C.~C. (2006). Spatially adaptive
  smoothing splines. {\em Biometrika\/} {\bf 93}, 113--125.

\bibitem[{Rue and Held(2005)}]{GMRFbook}
Rue, H. and Held, L. (2005). {\em Gaussian {M}arkov Random Fields: {T}heory and
  Applications\/}, volume 104 of {\em Monographs on Statistics and Applied
  Probability\/}. Chapman \& Hall, London.

\bibitem[{Rue et~al.(2009)Rue, Martino and Chopin}]{Rue:Mart:Chop:inla:2009}
Rue, H., Martino, S. and Chopin, N. (2009). Approximate {B}ayesian inference
  for latent {G}aussian models using integrated nested {L}aplace approximations
  (with discussion). {\em Journal of the Royal Statistical Society, Series B:
  Statistical Methodology\/} {\bf 71}, 319--392.

\bibitem[{Ruppert and Carroll(2000)}]{Rupp:Carr:spat:2000}
Ruppert, D. and Carroll, R.~J. (2000). Spatially-adaptive penalties for spline
  fitting. {\em Australian \& New Zealand Journal of Statistics\/} {\bf 42},
  205--223.

\bibitem[{Ruppert et~al.(2003)Ruppert, Wand and Carroll}]{rupp:wand:carr:2003}
Ruppert, D., Wand, M. and Carroll, R. (2003). {\em Semiparametric
  {R}egression\/}. Cambridge University Press, Cambridge.

\bibitem[{Scheipl and Kneib(2009)}]{scheipl:csda:09}
Scheipl, F. and Kneib, T. (2009). Locally adaptive {B}ayesian {P}-splines with
  a normal-exponential-gamma prior. {\em Computational Statistics and Data
  Analysis\/} {\bf 53}, 3533--3552.

\bibitem[{Silverman(1985)}]{silverman85}
Silverman, B.~W. (1985). Some aspects of the spline smoothing approach to
  non-parametric regression curve fitting. {\em Journal of the Royal
  Statistical Association Series B\/} {\bf 47}, 1--52.

\bibitem[{Simpson et~al.(2012)Simpson, Helton and Lindgren}]{simp:helt:lind:12}
Simpson, D., Helton, K. and Lindgren, F. (2012). On the connection between
  {O}'{S}ullivan splines, continuous random walk models, and smoothing splines.
  Technical report, Norwegian University of Science and Technology.

\bibitem[{Speckman and Sun(2003)}]{Spec:Sun:full:2003}
Speckman, P.~L. and Sun, D. (2003). Fully {B}ayesian spline smoothing and
  intrinsic autoregressive priors. {\em Biometrika\/} {\bf 90}, 289--302.

\bibitem[{Wahba(1978)}]{wahba78}
Wahba, G. (1978). Improper priors, spline smoothing and the problem of guarding
  against model errors in regression. {\em Journal of the Royal Statistical
  Society, Series B: Statistical Methodology\/} {\bf 40}, 364--372.

\bibitem[{Wahba(1990)}]{Wahb:spli:1990}
Wahba, G. (1990). {\em Spline Models for Observational Data\/}. SIAM [Society
  for Industrial and Applied Mathematics], Philadelphia.

\bibitem[{Walsh(1986)}]{walsh86}
Walsh, J. (1986). An introduction to stochastic partial differential equations.
  In {\em {\'E}cole d'{\'E}t{\'e} de Probabilit{\'e}s de Saint Flour XIV -
  1984\/} (R.~Carmona, H.~Kesten and J.~Walsh, eds.), volume 1180 of {\em
  Lecture Notes in Mathematics\/}, 265--439, Springer Berlin / Heidelberg,
  10.1007/BFb0074920.

\bibitem[{Wand and Ormerod(2008)}]{wand08}
Wand, M.~P. and Ormerod, J.~T. (2008). On semiparametric regression with
  {O}'{S}ullivan penalized splines. {\em Australian and New Zealand Journal of
  Statistics\/} {\bf 50}, 179--198.

\bibitem[{Wood et~al.(2002)Wood, Jiang and Tanner}]{Wood:Jiang:Tan:mix:2002}
Wood, S., Jiang, W. and Tanner, M. (2002). Bayesian mixture of splines for
  spatially adaptive nonparametric regression. {\em Biometrika\/} {\bf 89},
  513--528.

\bibitem[{Wood et~al.(2008)Wood, Kohn, Cottet, Jiang and Tanner}]{wood:kohn:08}
Wood, S.~A., Kohn, R., Cottet, R., Jiang, W. and Tanner, M. (2008). Locally
  adaptive nonparametric binary regression. {\em Journal of Computational and
  Graphical Statistics\/} {\bf 17}, 352--372.

\bibitem[{Yue et~al.(2012)Yue, Speckman and Sun}]{yue:spec:sun:aism:12}
Yue, Y., Speckman, P. and Sun, D. (2012). Priors for bayesian adaptive spline
  smoothing. {\em Annals of the Institute of Statistical Mathematics\/} {\bf
  64}, 577--613, 10.1007/s10463-010-0321-6.

\bibitem[{Yue and Speckman(2010)}]{yue:speck:non:10}
Yue, Y. and Speckman, P.~L. (2010). Nonstationary spatial {G}aussian {M}arkov
  random fields. {\em Journal of Computational and Graphical Statistics\/} {\bf
  19}, 96--116.

\bibitem[{Zhou and Shen(2001)}]{zhou01}
Zhou, S. and Shen, X. (2001). Spatially adaptive regression splines and
  accurate knot selection schemes. {\em Journal of the American Statistical
  Association\/} {\bf 96}, 247--259.

\end{thebibliography}
